\documentclass[11pt]{article}
\usepackage{hyperref,srcltx} 

\usepackage{amsmath}
\usepackage{amssymb}
\usepackage{amscd}
\usepackage{textcomp}

\usepackage{pstricks-add}
\usepackage{pstricks,pst-plot}

\usepackage{amsthm}
\theoremstyle{plain}

\theoremstyle{definition}

\theoremstyle{remark}

\newcounter{zahl}

\def\theenumi{(\alph{enumi})}

\def\p@enumii{\theenumi}

\long\def\forget#1{}

\newcommand{\DS}{\displaystyle}
\newcommand{\TS}{\textstyle}

\renewcommand{\phi}{\varphi}
\renewcommand{\epsilon}{\varepsilon}

\usepackage{amsfonts}
\newcommand{\BOne} {{\mathchoice{\hbox{\rm1\kern-2.7pt l\kern.9pt}}
                              {\hbox{\rm1\kern-2.7pt l\kern.9pt}}
                              {\hbox{\scriptsize\rm1\kern-2.3pt l\kern.4pt}}
                              {\hbox{\scriptsize\rm1\kern-2.4pt l\kern.5pt}}}}

\newcommand{\BD}{{\mathbb{D}}}

\newcommand{\BF}{{\mathbb{F}}}

\newcommand{\BQ}{{\mathbb{Q}}}
\newcommand{\BR}{{\mathbb{R}}}

\newcommand{\BZ}{{\mathbb{Z}}}

\newcommand{\es}{\enspace}
\newcommand{\open}{^\circ}

\newcommand{\ol}[1]{{\overline{#1}}}

\renewcommand{\angle}{\sphericalangle}

\usepackage{ifthen}

\newcommand{\invlim}[1][]{\ifthenelse{\equal{#1}{}}
{\DS \lim_{\longleftarrow}}
{\DS \lim_{\underset{#1}{\longleftarrow}}}
}

\newcommand{\dirlim}[1][]{\ifthenelse{\equal{#1}{}}
{\DS \lim_{\longrightarrow}}
{\DS \lim_{\underset{#1}{\longrightarrow}}}
}

\newcommand{\dotBD}{\vbox{\hbox{\kern2pt\bf.}\vskip-4.5pt\hbox{$\BD$}}}

\usepackage{color}
\newcounter{commentcounter}

\def\?{\ 
{\bf\color{red}???}\ 
\immediate\write16{}
\immediate\write16{Warning: There was still a question mark . . . }
\immediate\write16{}}

\newbox\mybox
\def\arrover#1{\mathrel{
       \setbox\mybox=\hbox spread 1.4em{\hfil$\scriptstyle#1$\hfil}
       \vbox{\offinterlineskip\copy\mybox
             \hbox to\wd\mybox{\rightarrowfill}}}}

\begin{document}

\author{Urs Hartl and Klaudia Kwickert}
\date{August 27, 2012}
\title{Constructing the Cubus simus and the Dodecaedron simum\\ via paper folding}

\maketitle

\begin{abstract}
The archimedean solids Cubus simus (snub cube) and Dodecaedron simum (snub dodecahedron) cannot be constructed by ruler and compass. We explain that for general reasons their vertices can be constructed via paper folding on the faces of a cube, respectively dodecahedron, and we present explicit folding constructions. The construction of the Cubus simus is particularly elegant. We also review and prove the construction rules of the other Archimedean solids.
\end{abstract}

\bigskip

\tableofcontents

%
%

\section{Construction of platonic and archimedean solids}

The \emph{archimedean solids} were first described by Archimedes in a now-lost work and were later rediscovered a few times. They are convex polyhedra with regular polygons as faces and many symmetries, such that their symmetry group acts transitively on their vertices \cite{AdamWyss, Cromwell}. They can be built by placing their vertices on the faces of a suitable \emph{platonic solid}, a \emph{tetrahedron}, a \emph{cube}, an \emph{octahedron}, a \emph{dodecahedron}, or an \emph{icosahedron}. In this article we discuss the constructibility of the archimedean solids by ruler and compass, or paper folding.

All the platonic solids can be constructed by ruler and compass. Such constructions are already explained by Euclid in Book XIII, Propositions~13--17, e.~g.~\cite{Euclid}. Of course the constructions have to be carried out in space. But in addition to the usual constructions with ruler and compass within a plane one only needs the following two extra constructions: to raise a perpendicular to a given plane through a point in this plane, and to construct a plane containing three (non-collinear) points.

The construction of the tetrahedron, the cube, and the octahedron in this manner is straight forward. To construct the dodecahedron and the icosahedron, one first needs to construct the regular pentagon or the golden ratio $\Phi=\frac{1+\sqrt{5}}{2} \approx 1.618$. Euclid even describes a construction of the tetrahedron together with the center of its circumsphere (Book XIII, Proposition 13).

\noindent\begin{minipage}[b]{0.68\linewidth}
Starting from the platonic solids, the archimedean solids can be obtained by constructing their vertices on the faces of the platonic solids. On the face of the tetrahedron the vertices of the \emph{truncated tetrahedron} are constructed by trisecting the edges. On the face of the cube the vertices of the \emph{truncated cube}, the \emph{cuboctahedron}, the \emph{rhombicuboctahedron}, and the \emph{truncated cuboctahedron} are constructed according to figure ~\ref{FigureCube}. 
\end{minipage}
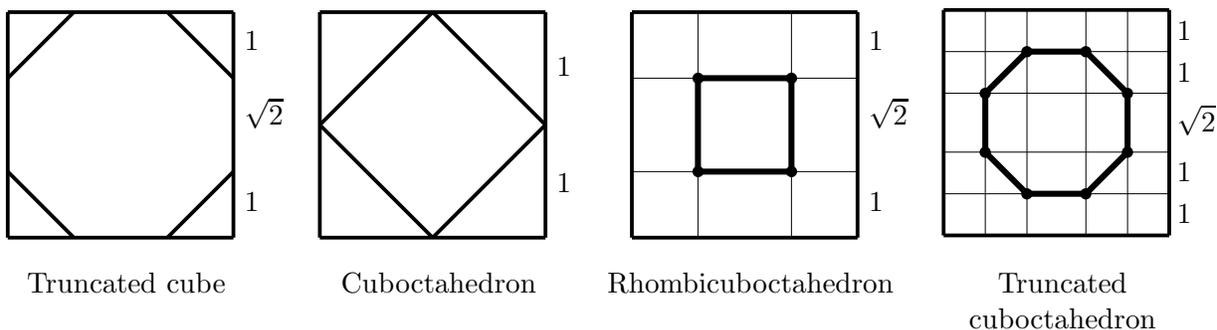
\begin{minipage}[b]{0.3\linewidth}
\begin{center}
\psset{unit=5mm}
\begin{pspicture}(0,6)(6,-0.5)
\psline[linewidth=0.5mm](3,5.196)(6,0)
\psline[linewidth=0.5mm](0,0)(3,5.196)
\psline[linewidth=0.5mm](0,0)(6,0)

\psline[linewidth=0.5mm](2,0)(1,1.732)
\psline[linewidth=0.5mm](4,0)(5,1.732)
\psline[linewidth=0.5mm](2,3.464)(4,3.464)

\put(5.7,1){$1$}
\put(4.9,2.6){$1$}
\put(4,4.2){$1$}
\end{pspicture}

Truncated tetrahedron
\end{center}
\end{minipage}

\begin{figure}[h]
\centering
\caption{Vertices of archimedean solids on the face of a cube.}
\label{FigureCube}

\psset{unit=0.5cm}

\begin{minipage}[t]{0.24\linewidth}
\begin{center}
\begin{pspicture}(0,6.8)(6.3,-0.5)
\psline[linewidth=0.5mm](0,0)(6,0)
\psline[linewidth=0.5mm](0,0)(0,6)
\psline[linewidth=0.5mm](0,6)(6,6)
\psline[linewidth=0.5mm](6,6)(6,0)

\psline[linewidth=0.5mm](1.757,0)(0,1.757)
\psline[linewidth=0.5mm](4.243,0)(6,1.757)
\psline[linewidth=0.5mm](0,4.243)(1.757,6)
\psline[linewidth=0.5mm](4.243,6)(6,4.243)

\put(6.3,0.7){$1$}
\put(6.3,3){$\sqrt{2}$}
\put(6.3,5){$1$}
\end{pspicture}

Truncated cube
\end{center}
\end{minipage}
\begin{minipage}[t]{0.24\linewidth}
\begin{center}
\begin{pspicture}(0,6.8)(6.3,-0.5)
\psline[linewidth=0.5mm](0,0)(6,0)
\psline[linewidth=0.5mm](0,0)(0,6)
\psline[linewidth=0.5mm](0,6)(6,6)
\psline[linewidth=0.5mm](6,6)(6,0)

\psline[linewidth=0.5mm](3,0)(6,3)
\psline[linewidth=0.5mm](6,3)(3,6)
\psline[linewidth=0.5mm](3,6)(0,3)
\psline[linewidth=0.5mm](0,3)(3,0)

\put(6.3,1.2){$1$}
\put(6.3,4.3){$1$}
\end{pspicture}

Cuboctahedron
\end{center}
\end{minipage}
\begin{minipage}[t]{0.24\linewidth}
\begin{center}
\psset{unit=0.5cm}
\begin{pspicture}(0,6.8)(6.3,-0.5)
\psline[linewidth=0.5mm](0,0)(6,0)
\psline[linewidth=0.5mm](0,0)(0,6)
\psline[linewidth=0.5mm](0,6)(6,6)
\psline[linewidth=0.5mm](6,6)(6,0)

\psline[linewidth=0.1mm](1.757,0)(1.757,6)
\psline[linewidth=0.1mm](4.243,0)(4.243,6)
\psline[linewidth=0.1mm](0,1.757)(6,1.757)
\psline[linewidth=0.1mm](0,4.243)(6,4.243)

\psline[linewidth=0.8mm](1.757,1.757)(1.757,4.243)
\psline[linewidth=0.8mm](1.757,4.243)(4.243,4.243)
\psline[linewidth=0.8mm](4.243,4.243)(4.243,1.757)
\psline[linewidth=0.8mm](4.243,1.757)(1.757,1.757)

\psdots[dotsize=1.5mm](1.757,1.757)(1.757,4.243)(4.243,4.243)(4.243,1.757)
\put(6.3,0.7){$1$}
\put(6.3,3){$\sqrt{2}$}
\put(6.3,5){$1$}
\end{pspicture}

Rhombicuboctahedron
\end{center}
\end{minipage}
%
%
\begin{minipage}[t]{0.24\linewidth}
\begin{center}
\psset{unit=0.554cm}
\begin{pspicture}(0,5.7)(5.7,-0.5)
\psline[linewidth=0.5mm](0,0)(5.414,0)
\psline[linewidth=0.5mm](5.414,0)(5.414,5.414)
\psline[linewidth=0.5mm](5.414,5.414)(0,5.414)
\psline[linewidth=0.5mm](0,5.414)(0,0)

\psline[linewidth=0.1mm](1,0)(1,5.414)
\psline[linewidth=0.1mm](2,0)(2,5.414)
\psline[linewidth=0.1mm](3.414,0)(3.414,5.414)
\psline[linewidth=0.1mm](4.414,0)(4.414,5.414)
\psline[linewidth=0.1mm](0,1)(5.414,1)
\psline[linewidth=0.1mm](0,2)(5.414,2)
\psline[linewidth=0.1mm](0,3.414)(5.414,3.414)
\psline[linewidth=0.1mm](0,4.414)(5.414,4.414)

\psline[linewidth=0.8mm](2,1)(3.414,1)
\psline[linewidth=0.8mm](3.414,1)(4.414,2)
\psline[linewidth=0.8mm](4.414,2)(4.414,3.414)
\psline[linewidth=0.8mm](4.414,3.414)(3.414,4.414)
\psline[linewidth=0.8mm](3.414,4.414)(2,4.414)
\psline[linewidth=0.8mm](2,4.414)(1,3.414)
\psline[linewidth=0.8mm](1,3.414)(1,2)
\psline[linewidth=0.8mm](1,2)(2,1)

\psdots[dotsize=1.5mm](2,1)(3.414,1)(4.414,2)(4.414,3.414)(3.414,4.414)(2,4.414)(1,3.414)(1,2)
\put(5.6,0.3){$1$}
\put(5.6,1.3){$1$}
\put(5.6,2.5){$\sqrt{2}$}
\put(5.6,3.7){$1$}
\put(5.6,4.7){$1$}
\end{pspicture}

Truncated cuboctahedron
\end{center}
\end{minipage}
\end{figure}
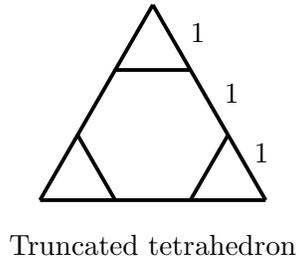

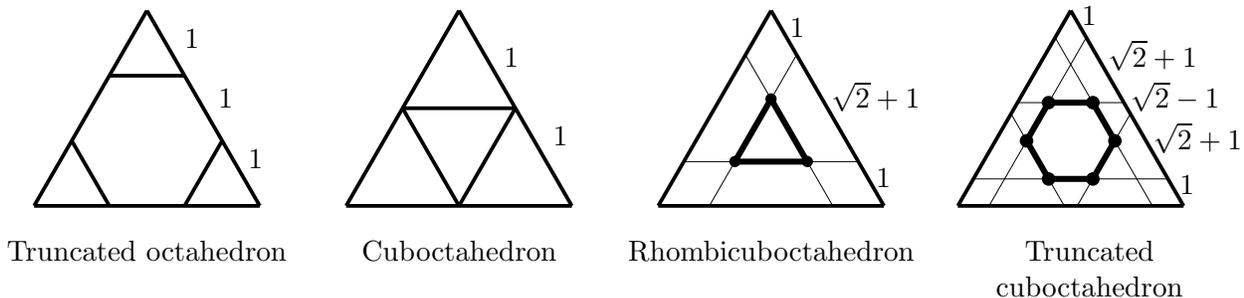
\begin{figure}[h]
\centering
\caption{Vertices of archimedean solids on the face of an octahedron.}
\label{FigureOctahedron}

\begin{minipage}[t]{0.24\linewidth}
\begin{center}
\psset{unit=5mm}
\begin{pspicture}(0,6)(6,-0.5)
\psline[linewidth=0.5mm](3,5.196)(6,0)
\psline[linewidth=0.5mm](0,0)(3,5.196)
\psline[linewidth=0.5mm](0,0)(6,0)

\psline[linewidth=0.5mm](2,0)(1,1.732)
\psline[linewidth=0.5mm](4,0)(5,1.732)
\psline[linewidth=0.5mm](2,3.464)(4,3.464)

\put(5.7,1){$1$}
\put(4.9,2.6){$1$}
\put(4,4.2){$1$}
\end{pspicture}

Truncated octahedron
\end{center}
\end{minipage}
\begin{minipage}[t]{0.24\linewidth}
\begin{center}
\psset{unit=5mm}
\begin{pspicture}(0,6)(6,-0.5)
\psline[linewidth=0.5mm](0,0)(3,5.196)
\psline[linewidth=0.5mm](0,0)(6,0)
\psline[linewidth=0.5mm](3,5.196)(6,0)

\psline[linewidth=0.5mm](3,0)(1.5,2.598)
\psline[linewidth=0.5mm](3,0)(4.5,2.598)
\psline[linewidth=0.5mm](1.5,2.598)(4.5,2.598)

\put(5.5,1.6){$1$}
\put(4.2,3.7){$1$}
\end{pspicture}

Cuboctahedron
\end{center}
\end{minipage}
\begin{minipage}[t]{0.24\linewidth}
\begin{center}
\psset{unit=5mm}
\begin{pspicture}(0,6)(6,-0.5)
\psline[linewidth=0.5mm](3,5.196)(6,0)
\psline[linewidth=0.5mm](0,0)(3,5.196)
\psline[linewidth=0.5mm](0,0)(6,0)

\psline[linewidth=0.1mm](1.359,0)(3.680,4.019)
\psline[linewidth=0.1mm](2.320,4.019)(4.641,0)
\psline[linewidth=0.1mm](0.680,1.177)(5.320,1.177)

\psline[linewidth=0.8mm](2.039,1.177)(3.961,1.177)
\psline[linewidth=0.8mm](3.961,1.177)(3,2.842)
\psline[linewidth=0.8mm](3,2.842)(2.039,1.177)

\psdots[dotsize=1.5mm](2.039,1.177)(3.961,1.177)(3,2.842)


\put(5.8,0.5){$1$}
\put(4.6,2.6){$\sqrt{2} + 1$}
\put(3.5,4.5){$1$}
\end{pspicture}

Rhombicuboctahedron
\end{center}
\end{minipage}
%
%
\begin{minipage}[t]{0.25\linewidth}
\begin{center}
\psset{unit=0.5cm}
\begin{pspicture}(0,6)(7,-0.5)
\psline[linewidth=0.5mm](3,5.196)(6,0)
\psline[linewidth=0.5mm](0,0)(3,5.196)
\psline[linewidth=0.5mm](0,0)(6,0)

\psline[linewidth=0.1mm](0.828,0)(3.414,4.479)
\psline[linewidth=0.1mm](2.586,4.479)(5.172,0)
\psline[linewidth=0.1mm](0.414,0.717)(5.586,0.717)

\psline[linewidth=0.1mm](3.172,0)(4.586,2.449)
\psline[linewidth=0.1mm](1.414,2.449)(2.828,0)
\psline[linewidth=0.1mm](1.586,2.747)(4.414,2.747)

\psline[linewidth=0.8mm](2.414,0.717)(3.586,0.717)
\psline[linewidth=0.8mm](3.586,0.717)(4.172,1.732)
\psline[linewidth=0.8mm](4.172,1.732)(3.586,2.747)
\psline[linewidth=0.8mm](3.586,2.747)(2.414,2.747)
\psline[linewidth=0.8mm](2.414,2.747)(1.828,1.732)
\psline[linewidth=0.8mm](1.828,1.732)(2.414,0.717)

\psdots[dotsize=1.8mm](2.414,0.717)(3.586,0.717)(4.172,1.732)(3.586,2.747)(2.414,2.747)(1.828,1.732)



\put(5.9,0.3){$1$}
\put(5.2,1.5){$\sqrt{2} + 1$}
\put(4.6,2.6){$\sqrt{2} - 1$}
\put(4,3.7){$\sqrt{2} + 1$}
\put(3.3,4.8){$1$}
\end{pspicture}

Truncated cuboctahedron
\end{center}
\end{minipage}
\end{figure}

On the face of the octahedron the vertices of the \emph{truncated octahedron}, the \emph{cuboctahedron}, the \emph{rhombicuboctahedron}, and the \emph{truncated cuboctahedron} are constructed according to figure ~\ref{FigureOctahedron}.

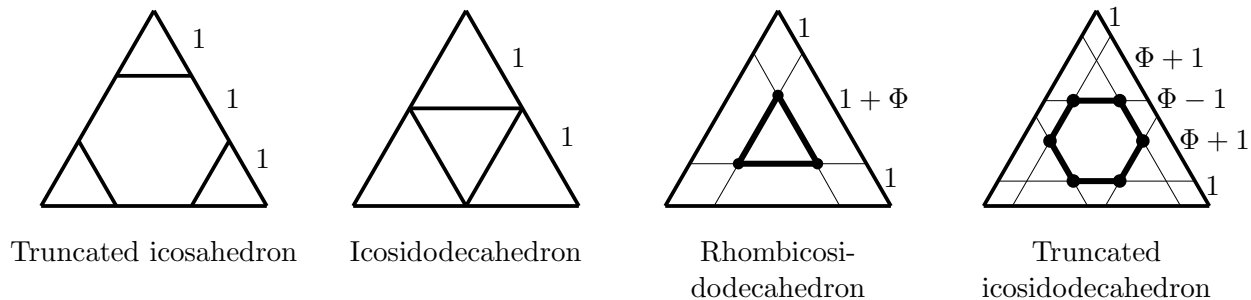
\begin{figure}[h]
\centering
\caption{Vertices of archimedean solids on the face of an icosahedron.}
\label{FigureIcosahedron}
\begin{minipage}[t]{0.24\linewidth}
\begin{center}
\psset{unit=5mm}
\begin{pspicture}(0,6)(6,-0.5)
\psline[linewidth=0.5mm](3,5.196)(6,0)
\psline[linewidth=0.5mm](0,0)(3,5.196)
\psline[linewidth=0.5mm](0,0)(6,0)

\psline[linewidth=0.5mm](2,0)(1,1.732)
\psline[linewidth=0.5mm](4,0)(5,1.732)
\psline[linewidth=0.5mm](2,3.464)(4,3.464)

\put(5.7,1){$1$}
\put(4.9,2.6){$1$}
\put(4,4.2){$1$}
\end{pspicture}

Truncated icosahedron
\end{center}
\end{minipage}
\begin{minipage}[t]{0.24\linewidth}
\begin{center}
\psset{unit=5mm}
\begin{pspicture}(0,6)(6,-0.5)
\psline[linewidth=0.5mm](0,0)(3,5.196)
\psline[linewidth=0.5mm](0,0)(6,0)
\psline[linewidth=0.5mm](3,5.196)(6,0)

\psline[linewidth=0.5mm](3,0)(1.5,2.598)
\psline[linewidth=0.5mm](3,0)(4.5,2.598)
\psline[linewidth=0.5mm](1.5,2.598)(4.5,2.598)

\put(5.5,1.6){$1$}
\put(4.2,3.7){$1$}
\end{pspicture}

Icosidodecahedron
\end{center}
\end{minipage}
\begin{minipage}[t]{0.24\linewidth}
\begin{center}
\psset{unit=5mm}
\begin{pspicture}(0,6)(6,-0.5)
\psline[linewidth=0.5mm](3,5.196)(6,0)
\psline[linewidth=0.5mm](0,0)(3,5.196)
\psline[linewidth=0.5mm](0,0)(6,0)

\psline[linewidth=0.1mm](1.299,0)(3.650,4.071)
\psline[linewidth=0.1mm](2.350,4.071)(4.641,0)
\psline[linewidth=0.1mm](0.650,1.125)(5.350,1.125)

\psline[linewidth=0.8mm](1.949,1.125)(4.051,1.125)
\psline[linewidth=0.8mm](4.051,1.125)(3,2.946)
\psline[linewidth=0.8mm](3,2.946)(1.949,1.125)

\psdots[dotsize=1.5mm](1.949,1.125)(4.051,1.125)(3,2.946)


\put(5.8,0.5){$1$}
\put(4.6,2.6){$1 + \Phi$}
\put(3.5,4.5){$1$}
\end{pspicture}

Rhombicosi-dodecahedron
\end{center}
\end{minipage}
%
%
\begin{minipage}[t]{0.25\linewidth}
\begin{center}
\psset{unit=0.5cm}
\begin{pspicture}(0,6)(6,-0.5)
\psline[linewidth=0.5mm](3,5.196)(6,0)
\psline[linewidth=0.5mm](0,0)(3,5.196)
\psline[linewidth=0.5mm](0,0)(6,0)

\psline[linewidth=0.1mm](0.764,0)(3.382,4.535)
\psline[linewidth=0.1mm](2.618,4.535)(5.236,0)
\psline[linewidth=0.1mm](0.382,0.662)(5.618,0.662)

\psline[linewidth=0.1mm](3.236,0)(4.618,2.394)
\psline[linewidth=0.1mm](1.382,2.394)(2.764,0)
\psline[linewidth=0.1mm](1.618,2.803)(4.382,2.803)

\psline[linewidth=0.8mm](2.382,0.662)(3.618,0.662)
\psline[linewidth=0.8mm](3.618,0.662)(4.236,1.732)
\psline[linewidth=0.8mm](4.236,1.732)(3.618,2.803)
\psline[linewidth=0.8mm](3.618,2.803)(2.382,2.803)
\psline[linewidth=0.8mm](2.382,2.803)(1.764,1.732)
\psline[linewidth=0.8mm](1.764,1.732)(2.382,0.662)

\psdots[dotsize=1.8mm](2.382,0.662)(3.618,0.662)(4.236,1.732)(3.618,2.803)(2.382,2.803)(1.764,1.732)



\put(5.9,0.3){$1$}
\put(5.2,1.5){$\Phi + 1$}
\put(4.6,2.6){$\Phi - 1$}
\put(4,3.7){$\Phi + 1$}
\put(3.3,4.8){$1$}
\end{pspicture}

Truncated icosidodecahedron
\end{center}
\end{minipage}
\end{figure}

On the face of the icosahedron the vertices of the \emph{truncated icosahedron}, the \emph{icosidodecahedron}, the \emph{rhombicosidodecahedron}, and the \emph{truncated icosidodecahedron} are constructed according to figure ~\ref{FigureIcosahedron}.

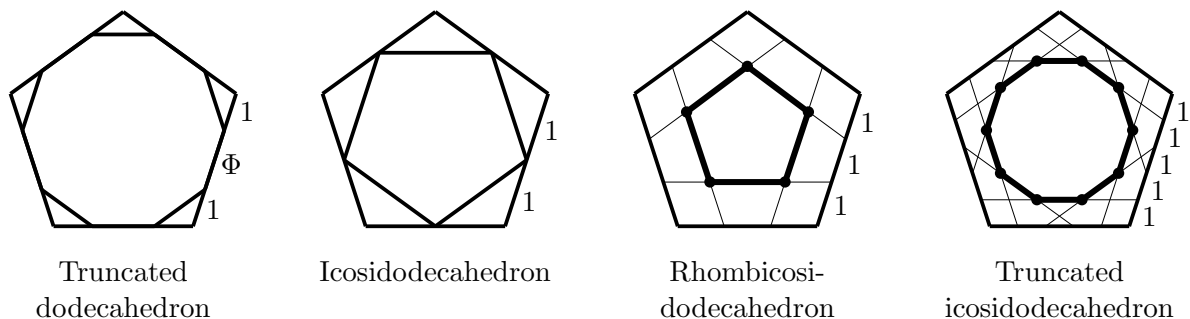
\begin{figure}[h]
\centering
\caption{Vertices of archimedean solids on the face of a dodecahedron.}
\label{FigureDodecahedron}
\begin{minipage}[t]{0.24\linewidth}
\begin{center}
\psset{unit=0.5cm}
\begin{pspicture}(0,6.5)(6,-0.5)
\psline[linewidth=0.5mm](1.146,0)(4.854,0)
\psline[linewidth=0.5mm](4.854,0)(6,3.527)
\psline[linewidth=0.5mm](6,3.527)(3,5.706)
\psline[linewidth=0.5mm](3,5.706)(0,3.527)
\psline[linewidth=0.5mm](0,3.527)(1.146,0)

\psline[linewidth=0.5mm](0.829,0.975)(2.171,0)
\psline[linewidth=0.5mm](2.171,0)(3.829,0)
\psline[linewidth=0.5mm](3.829,0)(5.171,0.975)
\psline[linewidth=0.5mm](5.171,0.975)(5.683,2.552)
\psline[linewidth=0.5mm](5.683,2.552)(5.171,4.129)
\psline[linewidth=0.5mm](5.171,4.129)(3.829,5.104)
\psline[linewidth=0.5mm](3.829,5.104)(2.171,5.104)
\psline[linewidth=0.5mm](2.171,5.104)(0.829,4.129)
\psline[linewidth=0.5mm](0.829,4.129)(0.317,2.552)
\psline[linewidth=0.5mm](0.317,2.552)(0.829,0.975)

\put(5.2,0.2){$1$}
\put(5.6,1.4){$\Phi$}
\put(6.1,2.8){$1$}
\end{pspicture}

Truncated dodecahedron
\end{center}
\end{minipage}
\begin{minipage}[t]{0.24\linewidth}
\begin{center}
\psset{unit=5mm}
\begin{pspicture}(0,6.5)(6,-0.5)
\psline[linewidth=0.5mm](1.146,0)(4.854,0)
\psline[linewidth=0.5mm](4.854,0)(6,3.527)
\psline[linewidth=0.5mm](6,3.527)(3,5.706)
\psline[linewidth=0.5mm](3,5.706)(0,3.527)
\psline[linewidth=0.5mm](0,3.527)(1.146,0)

\psline[linewidth=0.5mm](0.573,1.763)(3,0)
\psline[linewidth=0.5mm](3,0)(5.427,1.763)
\psline[linewidth=0.5mm](5.427,1.763)(4.5,4.617)
\psline[linewidth=0.5mm](4.5,4.617)(1.5,4.617)
\psline[linewidth=0.5mm](1.5,4.617)(0.573,1.763)

\put(5.3,0.4){$1$}
\put(5.9,2.3){$1$}
\end{pspicture}

Icosidodecahedron
\end{center}
\end{minipage}
\begin{minipage}[t]{0.24\linewidth}
\begin{center}
\psset{unit=0.5cm}
\begin{pspicture}(0,6.5)(6,-0.5)
\psline[linewidth=0.5mm](1.146,0)(4.854,0)
\psline[linewidth=0.5mm](4.854,0)(6,3.527)
\psline[linewidth=0.5mm](6,3.527)(3,5.706)
\psline[linewidth=0.5mm](3,5.706)(0,3.527)
\psline[linewidth=0.5mm](0,3.527)(1.146,0)

\psline[linewidth=0.1mm](2.382,0)(1,4.253)
\psline[linewidth=0.1mm](3.618,0)(5,4.253)
\psline[linewidth=0.1mm](0.382,2.315)(4,4.980)
\psline[linewidth=0.1mm](0.764,1.176)(5.236,1.176)
\psline[linewidth=0.1mm](5.618,2.315)(2,4.980)

\psline[linewidth=0.8mm](2,1.176)(4,1.176)
\psline[linewidth=0.8mm](4,1.176)(4.618,3.041)
\psline[linewidth=0.8mm](4.618,3.041)(3,4.254)
\psline[linewidth=0.8mm](3,4.254)(1.382,3.041)
\psline[linewidth=0.8mm](1.382,3.041)(2,1.176)

\psdots[dotsize=1.5mm](1.382,3.041)(2,1.176)(4,1.176)(4.618,3.041)(3,4.254)
\put(5.3,0.3){$1$}
\put(5.65,1.4){$1$}
\put(6,2.5){$1$}
\end{pspicture}

Rhombicosi-dodecahedron
\end{center}
\end{minipage}
\begin{minipage}[t]{0.24\linewidth}
\begin{center}
\psset{unit=0.5cm}
\begin{pspicture}(0,6.5)(6,-0.5)
\psline[linewidth=0.5mm](1.146,0)(4.854,0)
\psline[linewidth=0.5mm](4.854,0)(6,3.527)
\psline[linewidth=0.5mm](6,3.527)(3,5.706)
\psline[linewidth=0.5mm](3,5.706)(0,3.527)
\psline[linewidth=0.5mm](0,3.527)(1.146,0)

\psline[linewidth=0.1mm](1.888,0)(0.6,3.963)
\psline[linewidth=0.1mm](4.112,0)(5.4,3.963)
\psline[linewidth=0.1mm](0.229,2.822)(3.6,5.270)
\psline[linewidth=0.1mm](0.917,0.705)(5.083,0.705)
\psline[linewidth=0.1mm](5.771,2.822)(2.4,5.270)

\psline[linewidth=0.1mm](2.629,0)(5.542,2.116)
\psline[linewidth=0.1mm](3.371,0)(0.458,2.116)
\psline[linewidth=0.1mm](0.688,1.411)(1.8,4.834)
\psline[linewidth=0.1mm](1.2,4.399)(4.8,4.399)
\psline[linewidth=0.1mm](5.312,1.411)(4.2,4.834)

\psline[linewidth=0.8mm](2.4,0.705)(3.6,0.705)
\psline[linewidth=0.8mm](3.6,0.705)(4.570,1.410)
\psline[linewidth=0.8mm](4.570,1.410)(4.942,2.553)
\psline[linewidth=0.8mm](4.942,2.553)(4.571,3.693)
\psline[linewidth=0.8mm](4.571,3.693)(3.6,4.399)
\psline[linewidth=0.8mm](3.6,4.399)(2.4,4.399)
\psline[linewidth=0.8mm](2.4,4.399)(1.429,3.693)
\psline[linewidth=0.8mm](1.429,3.693)(1.058,2.553)
\psline[linewidth=0.8mm](1.058,2.553)(1.430,1.410)
\psline[linewidth=0.8mm](1.430,1.410)(2.4,0.705)

\psdots[dotsize=1.5mm](2.4,0.705)(3.6,0.705)(4.570,1.410)(4.942,2.553)(4.571,3.693)(2.4,4.399)(3.6,4.399)(1.429,3.693)(1.058,2.553)(1.430,1.410)
\put(5.2,0){$1$}
\put(5.425,0.7){$1$}
\put(5.65,1.4){$1$}
\put(5.875,2.1){$1$}
\put(6.1,2.8){$1$}
\end{pspicture}

Truncated icosidodecahedron
\end{center}
\end{minipage}
\end{figure}

Finally on the face of the dodecahedron the vertices of the \emph{truncated dodecahedron}, the \emph{icosidodecahedron}, the \emph{rhombicosidodecahedron}, and the \emph{truncated icosidodecahedron} are constructed according to figure ~\ref{FigureDodecahedron}. These construction rules are well known and easily verified (also see our appendix~\ref{AppConstructions}).

All these subdivisions of faces can be carried out with ruler and compass. In contrast the two remaining archimedean solids, \emph{Cubus simus} (also called \emph{snub cube}) and \emph{Dodecaedron simum} (also called \emph{snub dodecahedron}) cannot be constructed with ruler and compass only. This was proved by B.~Wei{\ss}bach and H.~Martini~\cite{WM}. They actually showed that the ratio $e:d$ between their edge length $e$ and the diameter $d$ of their circumsphere cannot be constructed by ruler and compass. Since one obtains this diameter and the edge length from the vertices, there can be no construction of the Cubus simus and the Dodecaedron simum with ruler and compass.

From an algebraic point of view the reason is that the inverse ratio $d:e$ satisfies the polynomial equation $(\frac{d}{e})^6 - 10(\frac{d}{e})^4 + 22(\frac{d}{e})^2 - 14 = 0$; see \cite[p.~126]{WM}. This equation is irreducible over $\BQ$ by Eisenstein's criterion (using the prime $2$) \cite[Chapter~11, Theorem~4.5]{Artin}, and hence $\frac{d}{e}$ generates a field extension $\BQ(\frac{d}{e})$ of degree $6$ over $\BQ$. On the other hand, with ruler and compass one can construct only ratios (or points with coordinates) lying in field extensions $K$ of $\BQ$ with $K\subset\BR$ whose degree is a power of $2$. More precisely, $K$ must contain subfields $K=K_n\supset K_{n-1}\supset\ldots\supset K_1\supset K_0=\BQ$ such that $K_i$ has degree $2$ over $K_{i-1}$ for all $i$; see \cite[Chapter~13, Theorem~4.9]{Artin}.

\medskip

Now there are various enhancements of construction techniques, allowing all constructions for which the degree of $K_i$ over $K_{i-1}$ is either $2$ or $3$. We want to focus on one of these enhancements namely \emph{paper folding}, see Martin~\cite[Chapter~10]{Martin}, in which one allows the following type of construction. Let two points $P$ and $Q$ and two lines $g$ and $h$ be drawn on a piece of paper but forbid the case where $P$ lies on $g$, $Q$ lies on $h$, and $g$ and $h$ are parallel or equal. We want to fold the paper creating one crease such that $P$ is placed onto $g$ and at the same time $Q$ is placed onto $h$. It turns out that there are at most three possibilities to achieve this and any crease obtained in this way will be a \emph{line constructible by paper folding}; see \cite[p.~149f]{Martin}. Mathematically the folding can be seen as a reflection at the crease. Notice that in the forbidden case, where $P$ already lies on $g$, $Q$ already lies on $h$, and $g$ and $h$ are parallel or equal, every fold along a crease orthogonal to $g$ and $h$ will achieve the goal. Since one does not want that all the lines orthogonal to a given line are constructible, this case is forbidden. All ruler and compass constructions can be carried out also by paper folding. Moreover, two of the classical Greek construction problems which are impossible by ruler and compass, the Delian problem to double the cube's volume, and the trisection of angles, can be solved by paper folding; see Martin~\cite[Lemmas~10.12 and 10.13 and Exercise 10.14]{Martin}.

In the next two sections we will describe paper folding methods to construct the vertices of the Cubus simus on the faces of a cube, and of the Dodecaedron simum on the faces of a dodecahedron. The construction of the Cubus simus is quite elegant, whereas our construction of the Dodecaedron simum is not. We are uncertain whether there is a nicer construction for the latter.

%
%

\section{The Cubus simus constructed out of a cube}\label{SectCubusSimus}

\begin{figure}[h]
\centering
\caption{The vertices of the Cubus simus on the faces of a cube (modified from \cite[p.~3]{Wyss}).}
\label{figure4.21}
\psset{unit=1.5cm}
\begin{pspicture}(-0.5,-0.5)(9,8)
\psline[linewidth=0.5mm](0,0)(6,0)
\psline[linewidth=0.5mm](0,0)(0,6)
\psline[linewidth=0.5mm](6,6)(6,0)
\psline[linewidth=0.5mm](0,6)(6,6)
\psline[linewidth=0.5mm](6,0)(9,1.5)
\psline[linewidth=0.5mm](9,1.5)(9,7.5)
\psline[linewidth=0.5mm](6,6)(9,7.5)
\psline[linewidth=0.5mm](3,7.5)(9,7.5)
\psline[linewidth=0.5mm](0,6)(3,7.5)

\psline[linestyle=dotted](2.113,0)(2.113,1.369)
\psline[linestyle=dotted](2.113,1.369)(0,1.369)
\psline[linestyle=dotted](4.631,0)(4.631,2.113)
\psline[linestyle=dotted](4.631,2.113)(6,2.113)
\psline[linestyle=dotted](0,3.887)(1.369,3.887)
\psline[linestyle=dotted](1.369,3.887)(1.369,6)
\psline[linestyle=dotted](6,4.631)(3.887,4.631)
\psline[linestyle=dotted](3.887,4.631)(3.887,6)
\psline[linestyle=dotted](6,3.887)(6.684,4.229)
\psline[linestyle=dotted](6.684,4.228)(6.684,6.342)
\psline[linewidth=0.8mm](3.887,4.631)(4.631,2.113)
\psline[linewidth=0.8mm](2.113,1.369)(4.631,2.113)
\psline[linewidth=0.8mm](2.113,1.369)(1.369,3.887)
\psline[linewidth=0.8mm](1.369,3.887)(3.887,4.631)
\psline(6.684,4.228)(3.887,4.631)
\psline(6.684,4.228)(4.631,2.113)
\psline(4.631,2.113)(7.057,1.897)
\psline[linewidth=0.8mm](6.684,4.228)(7.057,1.897)
\psline[linewidth=0.8mm](7.057,1.897)(8.316,3.271)
\psline[linewidth=0.8mm](8.316,3.271)(7.943,5.603)
\psline[linewidth=0.8mm](7.943,5.603)(6.684,4.228)
\psline(1.369,3.887)(2.798,6.342)
\psline(2.798,6.342)(3.887,4.631)
\psline[linewidth=0.8mm](2.798,6.342)(5.688,6.528)
\psline[linewidth=0.8mm](5.688,6.528)(6.202,7.158)
\psline[linewidth=0.8mm](6.202,7.158)(3.312,6.972)
\psline[linewidth=0.8mm](3.312,6.972)(2.798,6.342)
\psline[linestyle=dotted](5.688,6.528)(4.631,6)
\psline[linestyle=dotted](5.688,6.528)(7.057,6.528)
\psline(5.688,6.528)(7.943,5.603)
\psline(5.688,6.528)(6.684,4.228)
\psline(5.688,6.528)(3.887,4.631)
\put(3.6,4.2){$A_1$}
\put(4.1,2.2){$A_2$}
\put(2.1,1.6){$A_3$}
\put(1.5,3.7){$A_4$}
\put(6.8,4){$B_1$}
\put(7.9,5.7){$B_2$}
\put(8.4,3.1){$B_3$}
\put(7.1,1.7){$B_4$}
\put(5.4,6.7){$C_1$}
\put(2.4,6.3){$C_2$}
\put(3,7.1){$C_3$}
\put(6.3,7.2){$C_4$}
\put(1.1,5){$x$}
\put(0.5,3.7){$y$}
\put(1,1.2){$x$}
\put(2.2,0.6){$y$}
\put(4.7,0.8){$x$}
\put(5.5,2.3){$y$}
\put(5.3,4.7){$x$}
\put(4,5.5){$y$}
\put(6.4,6.65){$y$}
\put(4.6,6.1){$x$}
\put(6.45,5.2){$x$}
\put(6.1,4.1){$y$}
\end{pspicture}
\end{figure}
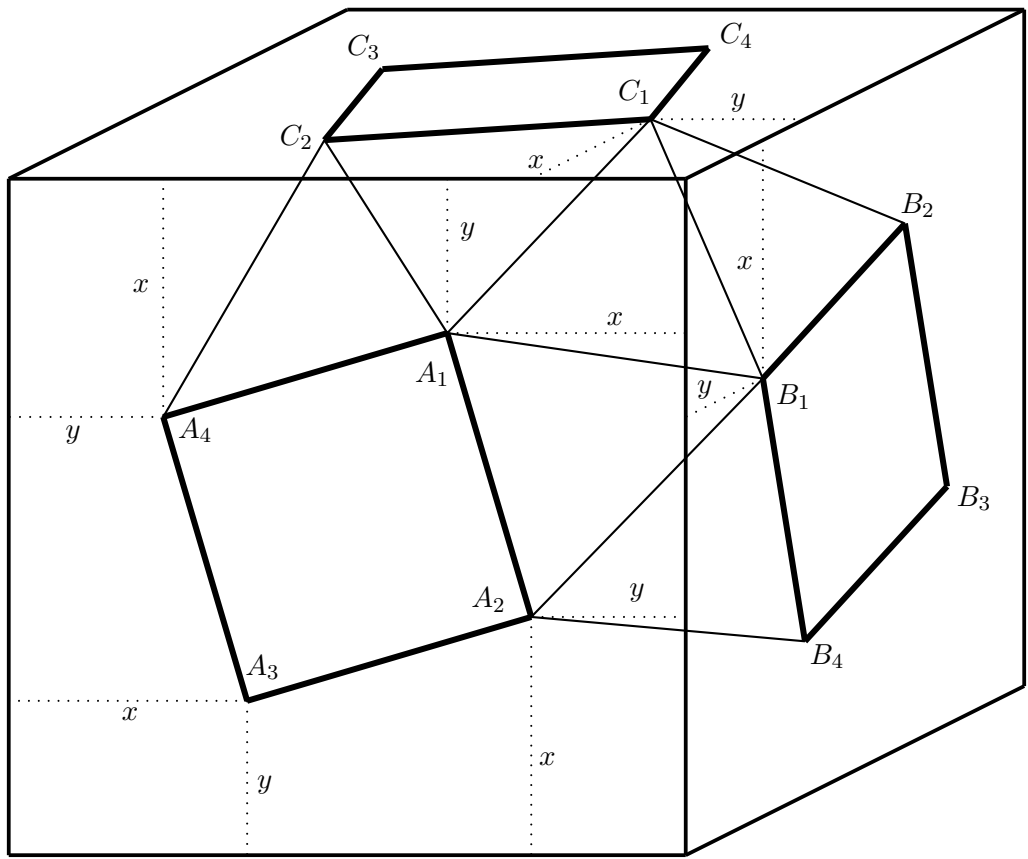

The \emph{Cubus simus} (or \emph{snub cube}) is constructed by placing smaller squares on the faces of the cube which are slightly rotated either to the left or to the right. Then the vertices of the smaller rotated squares are connected such that each vertex belongs to one square and four triangles; see figure~\ref{figure4.21}.

We normalize the edge length of the cube to $1$. Let $x, y$ be the distances of a vertex $A_1$ from the edges of the cube (i.e.~the lengths of the perpendiculars from $A_1$ onto these edges) such that $x \ge y$. By reasons of symmetry the distances $x, y$ must be the same for all the vertices of the Cubus simus. From Pythagoras theorem we compute
\begin{eqnarray*}
|A_1 A_2|^2 & = & (x-y)^2 + (1-x-y)^2 \\
|A_1 B_1|^2 & = & (x-y)^2 + x^2 + y^2 \\
|A_2 B_1|^2 & = & (1-2x)^2 + y^2 + y^2 \ .
\end{eqnarray*}
Since the triangle $\Delta A_1 A_2 B_1$ should be equilateral we obtain the two equations
\begin{equation*}
|A_1 A_2|^2 = |A_1 B_1|^2 \quad \textrm{and} \quad |A_1 A_2|^2 = |A_2 B_1|^2 \ .
\end{equation*}
The first implies $(x-y)^2 + x^2 + y^2 = (x-y)^2 + (1-x-y)^2 = (x-y)^2 + 1 + x^2 + y^2 - 2x - 2y + 2xy$, whence
\begin{equation}\label{eq:410}
1 - 2x - 2y + 2xy = 0 \ .
\end{equation}
The second implies
\begin{equation*}
x^2 - 2xy + y^2 + 1 + x^2 + y^2 - 2x - 2y + 2xy = 1 - 4x + 4x^2 + 2y^2,
\end{equation*}
whence $2x^2 + 2y^2 + 1 - 2x - 2y = 1 - 4x + 4x^2 + 2y^2$, and $-2x^2 + 2x = 2y$, and
\begin{equation}\label{eq:411}
y = x - x^2 \ .
\end{equation}
We substitute \eqref{eq:411} into \eqref{eq:410} to obtain 
\begin{equation}\label{eq:412}
2x^3 - 4x^2 + 4x - 1 = 0 \ .
\end{equation}
By Cardano's formula for the solutions of degree three equations we find
\begin{eqnarray*}
x & = & \frac{2}{3}+\frac{1}{6} \cdot \left ( \sqrt[3]{- 26 + 6 \cdot \sqrt[2]{33}} - \sqrt[3]{26 + 6 \cdot \sqrt[2]{33}} \right ) \\
& \approx & 0,352201129 \\
y & \approx & 0,228155494 \ .
\end{eqnarray*}

From equation~\eqref{eq:412} one obtains another proof for the fact that the vertices of the Cubus simus cannot be constructed with ruler and compass on the face of the cube. Namely \eqref{eq:412} is irreducible in $\BQ[x]$. Indeed, if it were reducible then it would have a solution $x\in\BQ$ which we write as a reduced fraction $x=\frac{p}{q}$ with relatively prime integers $p,q$. Multiplying with $q^3$ yields
\[
2p^3 - 4qp^2 + 4 q^2p - q^3 = 0
\]
and one sees that $q$ must be divisible by $2$, say $q=2r$. Substituting this and dividing by $2$ gives
\[
p^3 - 4rp^2 + 8r^2p-4r^3 =0\,.
\]
It follows that also $p$ is divisible by $2$. But this contradicts that $p$ and $q$ are relatively prime. Thus \eqref{eq:412} is irreducible and $\BQ(x)$ is a field extension of degree $3$ over $\BQ$. This proves that $x$ cannot be constructed by ruler and compass.

\bigskip

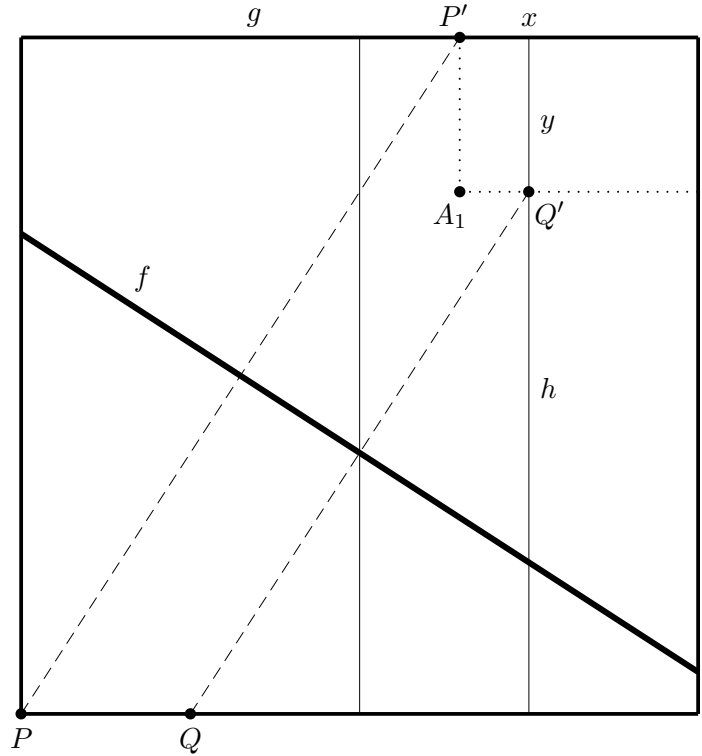
\begin{figure}[h]
\centering
\caption{Constructing the vertices of the Cubus simus by paper folding.}
\label{FigureFoldingCubusSimus}
\psset{unit=1.5cm}
\begin{pspicture}(-0.5,-0.5)(6,6.8)
\psline[linewidth=0.5mm](0,0)(6,0)
\psline[linewidth=0.5mm](0,0)(0,6)
\psline[linewidth=0.5mm](6,6)(6,0)
\psline[linewidth=0.5mm](0,6)(6,6)
\psline[linewidth=0.1mm](3,0)(3,6)
\psline[linewidth=0.1mm](4.5,0)(4.5,6)
\psline[linewidth=0.8mm](0,4.259)(6,0.372)
\psline[linestyle=dashed, linewidth=0.1mm](0,0)(3.887,6)
\psline[linestyle=dashed, linewidth=0.1mm](1.5,0)(4.5,4.631)
\psline[linestyle=dotted](6,4.631)(3.887,4.631)
\psline[linestyle=dotted](3.887,4.631)(3.887,6)
\psdots[dotsize=1.5mm](0,0)(1.5,0)(3.887,6)(4.5,4.631)(3.887,4.631)
\put(-0.1,-0.3){$P$}
\put(1.4,-0.3){$Q$}
\put(3.7,6.1){$P'$}
\put(4.55,4.35){$Q'$}
\put(3.65,4.35){$A_1$}
\put(1,3.8){$f$}
\put(2,6.15){$g$}
\put(4.6,2.8){$h$}
\put(4.43,6.1){$x$}
\put(4.6,5.2){$y$}
\end{pspicture}
\end{figure}

However it can be constructed by paper folding. Begin with a square face of the cube and divide one of its edges into four equal segments. Let $P$ be one endpoint of this edge and let $Q$ be the first division point. Let $g$ be the opposite edge of the square. Construct $h$ as the perpendicular of the last division point onto $g$ as in figure~\ref{FigureFoldingCubusSimus}. Then fold the square such that $P$ comes to lie on $g$ and $Q$ comes to lie on $h$. In figure~\ref{FigureFoldingCubusSimus} the crease $f$ is marked as a bold line and $P$ and $Q$ come to lie on $P'$, respectively $Q'$. We claim that the distance of $P'$ from the nearest corner of the cube's face is $x$ and the distance of $Q'$ from $g$ is $y$.

Before we prove the claim, notice that in this construction no ruler and compass are needed. The division of the edge into quarters and the construction of the line $h$ can be achieved by three times folding in halves. Finally $A_1$ is obtained as the image of $Q'$ under the folding which places $P'$ onto the line $h$ and the line $g$ onto itself.

To prove the claim we fix coordinates on the square such that $P={0\choose 0}, Q={1/4\choose 0}$, and the lines are given as $g=\{{0\choose 1}+s\cdot{1\choose 0}:s\in\BR\}$ and $h=\{{3/4\choose 0}+t\cdot{0\choose 1}:t\in\BR\}$. For the fold $f$ we set $f=\{{0\choose b}+u\cdot{1\choose c}:u\in\BR\}$ with unknown $b,c$. Folding along $f$ can be viewed as the reflection at $f$ which in our coordinates is given as the map
\[
R_f:{v_1\choose v_2}\longmapsto{0\choose b}+{\textstyle\frac{1}{1+c^2}}\left(\begin{array}{cc} 1-c^2&2c\\2c&c^2-1 \end{array}\right)\cdot\bigl({v_1\choose v_2}-{0\choose b}\bigr)\,.
\]
This reflection sends $P$ and $Q$ to 
\begin{eqnarray*}
P'\es=\es R_f(P)&=&R_f({0\choose 0})\es =\es {0\choose b}+{\textstyle\frac{1}{1+c^2}}{-2bc\choose -bc^2+b}\qquad\text{and}\\[3mm]
Q'\es=\es R_f(Q)&=&R_f({1/4\choose 0})\es =\es {0\choose b}+{\textstyle\frac{1}{1+c^2}}{-2bc+1/4-c^2/4\choose -bc^2+b+c/2}\,.
\end{eqnarray*}
Our requirements that $P'$ lies on $g$ and $Q'$ lies on $h$ mean $P'={1-d\choose 1}$ and $Q'={3/4\choose 1-e}$ for some $d,e$. This translates into
\begin{eqnarray*}
\textstyle b+\frac{1}{1+c^2}(-bc^2+b)&=&1\qquad\text{and}\\[2mm]
\textstyle \frac{1}{1+c^2}(-2bc+1/4-c^2/4)&=&3/4\,,
\end{eqnarray*}
whence
\begin{eqnarray}
\label{EqForbc1} 2b&=&1+c^2\qquad\text{and}\\[2mm]
\label{EqForbc2} -2bc&=&c^2+1/2\,.
\end{eqnarray}
Plugging \eqref{EqForbc1} into \eqref{EqForbc2} we obtain $c^3+c^2+c+1/2=0$. Writing this as an equation for $c+1$ yields $(c+1)^3-2(c+1)^2+2(c+1)-1/2=0$, which is the same as our old equation \eqref{eq:412} for $x$. Since the derivative of \eqref{eq:412} is $6x^2-8x+4=2x^2+4(x-1)^2$ which is always positive for real numbers $x$, we conclude that \eqref{eq:412} has exactly one real solution. Therefore $c+1=x$ and $c=x-1=\frac{y}{x}$. Using $1+c^2=x^2-2x+2=\frac{1}{2x}$ we can now compute the first coordinate of $P'={1-d\choose 1}$ as follows.
\[
d=1-\frac{-2bc}{1+c^2}=\frac{1+c^2-(c^2+1/2)}{1+c^2}=\frac{1}{2(1+c^2)}=x\,.
\]
Also we compute the second coordinate of $Q'={3/4\choose 1-e}$ as follows.
\[
e=1-b-\frac{-bc^2+b+c/2}{1+c^2}=\frac{1+c^2-b-bc^2+bc^2-b-c/2}{1+c^2}=\frac{-c}{2(1+c^2)}=x\cdot\frac{y}{x}=y\,.
\]
Altogether this proves our claim that the distance of $P'$ from the nearest corner ${1\choose 1}$ of the cube's face is $x$ and the distance of $Q'$ from $g$ is $y$.\qed

%
%

\section{The Dodecaedron simum constructed out of a dodecahedron}

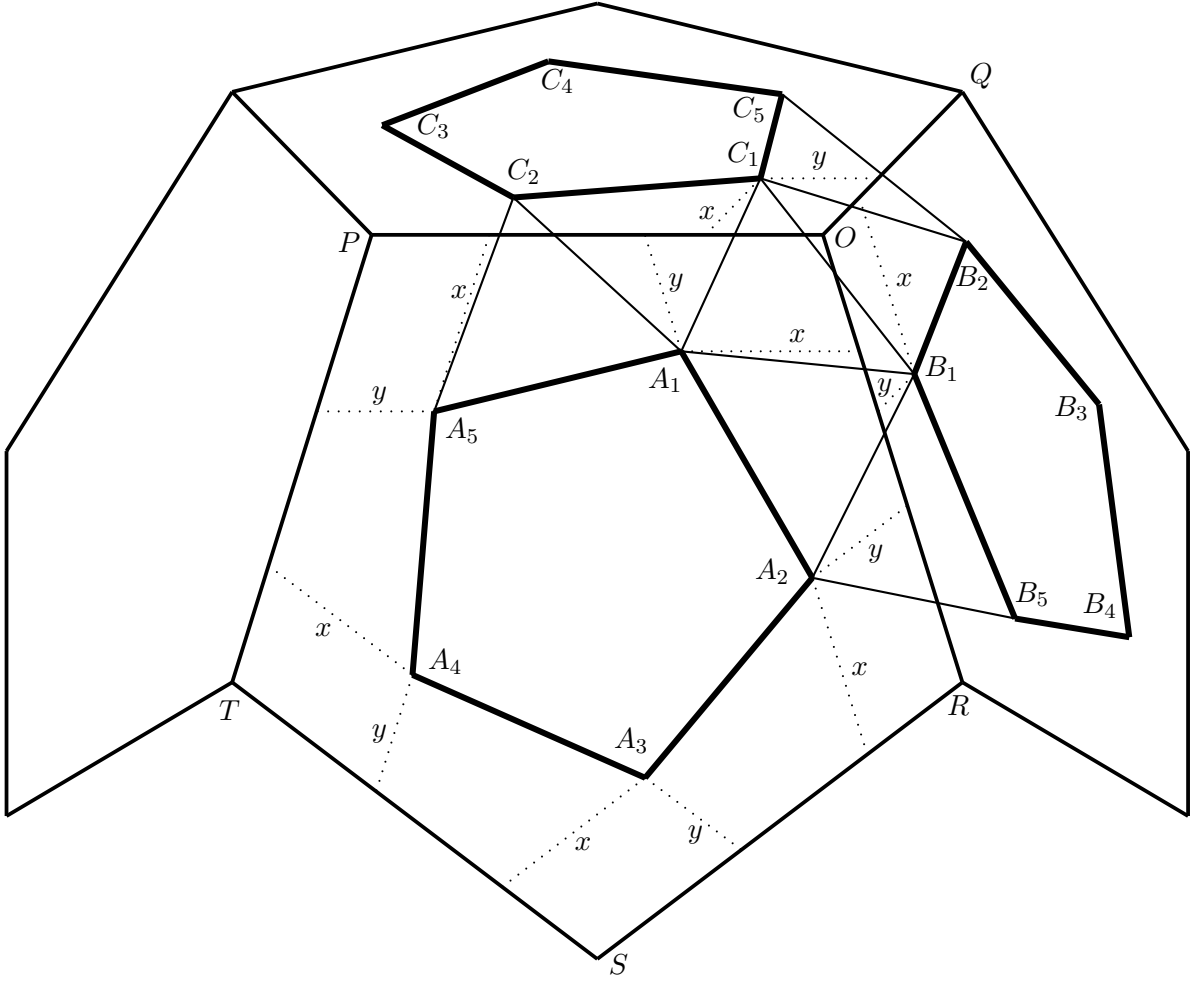
\begin{figure}[h]
\centering
\caption{The vertices of the Dodecaedron simum on the faces of a dodecahedron.}
\label{figure4.23}
\psset{unit=1.5cm}
\begin{pspicture}(-2.5,-0.5)(8.5,8.5)
\psline[linewidth=0.5mm](5, 6.155)(6.236, 2.186)
\psline[linewidth=0.5mm](6.236, 2.186)(3,-0.268)
\psline[linewidth=0.5mm](3,-0.268)(-0.236, 2.186)
\psline[linewidth=0.5mm](-0.236, 2.186)(1, 6.155)
\psline[linewidth=0.5mm](1, 6.155)(5, 6.155)
\put(5.1,6.05){$O$}
\put(0.7,6.0){$P$}
\put(6.3,7.5){$Q$}
\put(6.1,1.9){$R$}
\put(3.1,-0.4){$S$}
\put(-0.35, 1.85){$T$}

\psline[linewidth=0.5mm](6.236,2.186)(8.236,1)
\psline[linewidth=0.5mm](8.236,1)(8.236,4.237)
\psline[linewidth=0.5mm](8.236,4.237)(6.236,7.423)
\psline[linewidth=0.5mm](6.236,7.423)(5,6.155)
\psline[linewidth=0.5mm](6.236,7.423)(3,8.207)
\psline[linewidth=0.5mm](3,8.207)(-0.236,7.423)
\psline[linewidth=0.5mm](-0.236,7.423)(1,6.155)
\psline[linewidth=0.5mm](-0.236,2.186)(-2.236,1)
\psline[linewidth=0.5mm](-2.236,1)(-2.236,4.237)
\psline[linewidth=0.5mm](-2.236,4.237)(-0.236,7.423)
\psline[linewidth=0.8mm](3.744,5.122)(1.554,4.589)
\psline[linewidth=0.8mm](1.554,4.589)(1.362,2.252)
\psline[linewidth=0.8mm](1.362,2.252)(3.424,1.339)
\psline[linewidth=0.8mm](3.424,1.339)(4.906,3.113)
\psline[linewidth=0.8mm](4.906,3.113)(3.744,5.122)
\put(3.45,4.8){$A_1$}
\put(4.4,3.1){$A_2$}
\put(3.15,1.6){$A_3$}
\put(1.5,2.3){$A_4$}
\put(1.65,4.35){$A_5$}
\psline[linestyle=dotted](3.422,6.155)(3.744,5.122)
\put(3.63,5.7){$y$}
\psline[linestyle=dotted](3.744,5.122)(5.322,5.122)
\put(4.7,5.2){$x$}
\psline[linestyle=dotted](0.512,4.589)(1.554,4.589)
\put(1,4.7){$y$}
\psline[linestyle=dotted](1.554,4.589)(2.042,6.155)
\put(1.7,5.6){$x$}
\psline[linestyle=dotted](1.040,1.218)(1.362,2.252)
\put(1,1.7){$y$}
\psline[linestyle=dotted](1.362,2.252)(0.086,3.220)
\put(0.5,2.6){$x$}
\psline[linestyle=dotted](4.277,0.700)(3.424,1.339)
\put(3.8,0.8){$y$}
\psline[linestyle=dotted](3.424,1.339)(2.157,0.371)
\put(2.8,0.7){$x$}
\psline[linestyle=dotted](5.748,3.752)(4.906,3.113)
\put(5.4,3.3){$y$}
\psline[linestyle=dotted](4.906,3.113)(5.393,1.546)
\put(5.25,2.25){$x$}
\psline[linewidth=0.8mm](5.809,4.920)(6.703,2.752)
\psline[linewidth=0.8mm](6.703,2.752)(7.715,2.586)
\psline[linewidth=0.8mm](7.715,2.586)(7.447,4.651)
\psline[linewidth=0.8mm](7.447,4.651)(6.269,6.093)
\psline[linewidth=0.8mm](6.269,6.093)(5.809,4.920)
\put(5.9,4.9){$B_1$}
\put(6.17,5.7){$B_2$}
\put(7.05,4.55){$B_3$}
\put(7.3,2.8){$B_4$}
\put(6.7,2.9){$B_5$}
\psline[linestyle=dotted](5.488,4.589)(5.809,4.920)
\put(5.48,4.77){$y$}
\psline[linestyle=dotted](5.809,4.920)(5.322,6.486)
\put(5.65,5.7){$x$}
\psline[linewidth=0.8mm](4.446,6.656)(4.638,7.402)
\psline[linewidth=0.8mm](4.638,7.402)(2.566,7.694)
\psline[linewidth=0.8mm](2.566,7.694)(1.094,7.127)
\psline[linewidth=0.8mm](1.094,7.127)(2.256,6.486)
\psline[linewidth=0.8mm](2.256,6.486)(4.446,6.656)
\put(4.15,6.8){$C_1$}
\put(2.2,6.65){$C_2$}
\put(1.4,7.05){$C_3$}
\put(2.5,7.45){$C_4$}
\put(4.2,7.2){$C_5$}
\psline[linestyle=dotted](5.488,6.656)(4.446,6.656)
\put(4.9,6.78){$y$}
\psline[linestyle=dotted](4.446,6.656)(3.958,6.155)
\put(3.9,6.3){$x$}
\psline(3.744,5.122)(5.809,4.920)
\psline(5.809,4.920)(4.446,6.656)
\psline(4.446,6.656)(3.744,5.122)
\psline(3.744,5.122)(2.256,6.486)
\psline(2.256,6.486)(1.554,4.589)
\psline(5.809,4.920)(4.906,3.113)
\psline(4.906,3.113)(6.703,2.752)
\psline(4.446,6.656)(6.269,6.093)
\psline(6.269,6.093)(4.638,7.402)
\end{pspicture}
\end{figure}

The \emph{Dodecaedron simum} (or \emph{snub dodecahedron}) is constructed by placing smaller pentagons on the faces of a dodecahedron which are slightly rotated either to the left or to the right. Then the vertices of the smaller rotated pentagons are connected such that each vertex belongs to one pentagon and four triangles; see figure~\ref{figure4.23}.

We normalize the edge length of the dodecahedron to $1$. Then the vectors $\overrightarrow{OP}$, $\overrightarrow{OQ}$, and $\overrightarrow{OR}$ all have length $1$. Since all three angles at $O$ have the same size $\angle POQ=\angle POR =\angle QOR=108\open$ we compute for the scalar product
\begin{equation}\label{EqScalarProd}
\langle\overrightarrow{OP},\overrightarrow{OQ}\rangle\es=\es \langle\overrightarrow{OP},\overrightarrow{OR}\rangle\es=\es \langle\overrightarrow{OQ},\overrightarrow{OR}\rangle\es=\es \cos 108\open\es=\es -\frac{1}{2\Phi}\,,
\end{equation}
where again $\Phi=\frac{1+\sqrt{5}}{2} \approx 1.618$ is the golden ratio.

Join the vertex $A_1$ of the Dodecaedron simum by line segments with the edges $[OP]$ and $[OR]$ of the dodecahedron to form a parallelogram with the vertex $O$ (the dotted lines labeled $x$ and $y$ in figure~\ref{figure4.23}). Let $x, y$ be the lengths of these segments such that $x \ge y$. By reasons of symmetry the distances $x, y$ must be the same for all the vertices of the Dodecaedron simum. Since the ``diagonal'' $\overrightarrow{PS}$ of the dodecahedron's face has length $\Phi$ and is parallel to $\overrightarrow{OR}$ we have
\[
\overrightarrow{OS}\es=\es\overrightarrow{OP}+\Phi\cdot\overrightarrow{OR}\qquad\text{and}\qquad\overrightarrow{RS}\es=\es\overrightarrow{OS}-\overrightarrow{OR}\es=\es\overrightarrow{OP}+(\Phi-1)\cdot\overrightarrow{OR}\,.\\[2mm]
\]
Then the vertices of the Dodecaedron simum have the following coordinates
\[
\begin{array}{lclcl}
\overrightarrow{OA_1}&=&x\cdot\overrightarrow{OP}+y\cdot\overrightarrow{OR}\,,\\[2mm]
\overrightarrow{OA_2}&=&\overrightarrow{OR}-x\cdot\overrightarrow{OR}+y\cdot\overrightarrow{RS}&=&y\cdot\overrightarrow{OP}+(1-x-y+y\Phi)\cdot\overrightarrow{OR}\,,\\[2mm]
\overrightarrow{OB_1}&=&x\cdot\overrightarrow{OR}+y\cdot\overrightarrow{OQ}\,.\\[2mm]
\end{array}
\]
Using the values for the scalar product from \eqref{EqScalarProd} we compute the distances 
\begin{eqnarray*}
|\overrightarrow{A_1A_2}|^2&=&\bigl|(y-x)\cdot\overrightarrow{OP}+(1-x-2y+y\Phi)\cdot\overrightarrow{OR}\bigr|^2\\[2mm]
&=&(3-\Phi)x^2+(3-\Phi)y^2+(4-3\Phi)xy-(3-\Phi)x-(3-\Phi)y+1\,,\\[2mm]
|\overrightarrow{A_1B_1}|^2&=&\bigl|-x\cdot\overrightarrow{OP}+(x-y)\cdot\overrightarrow{OR}+y\cdot\overrightarrow{OQ}\bigr|^2\\[2mm]
&=&(1+\Phi)(x^2+y^2-xy)\,,\\[2mm]
|\overrightarrow{B_1A_2}|^2&=&\bigl|y\cdot\overrightarrow{OP}+(1-2x-y+y\Phi)\cdot\overrightarrow{OR}-y\cdot\overrightarrow{OQ}\bigr|^2\\[2mm]
&=&4x^2+3y^2+4(1-\Phi)xy-4x+2(\Phi-1)y+1\,.
\end{eqnarray*}
Since the triangle $\Delta A_1 A_2 B_1$ should be equilateral we obtain the two equations
\begin{equation*}
|A_1 A_2|^2 = |A_2 B_1|^2  \quad \textrm{and} \quad |A_1 B_1|^2 = |A_2 B_1|^2
\end{equation*}
which give us
\begin{eqnarray}
\label{Eq1ForDodSimum}y^2&=&-\Phi x^2+\Phi x-\Phi y+xy\quad\text{and}\\[2mm]
\label{Eq2ForDodSimum}(\Phi-2)y^2&=&(3-\Phi)x^2-4x+1+(5-3\Phi)xy+2(\Phi-1)y
\end{eqnarray}
We substitute \eqref{Eq1ForDodSimum} into \eqref{Eq2ForDodSimum} to obtain 
\begin{equation}\label{EqYForDodSimum}
y\es=\es\frac{(4-2\Phi)x^2+(\Phi-5)x+1}{(4\Phi-7)x+1-\Phi}\,.
\end{equation}
Substituting this back into \eqref{Eq2ForDodSimum} yields
\[
\bigl((25\Phi-41)x^3+(30-17\Phi)x^2-5x+\Phi\bigr)\cdot x\es=\es 0\,.
\]
In the case $x=0$, which implies $y=-\Phi$, the point $A_1$ does not lie on the face of the dodecahedron. Therefore we must have $x\ne0$ and
\begin{equation}\label{EqXForDodSimum}
(25\Phi-41)x^3+(30-17\Phi)x^2-5x+\Phi\es=\es 0
\end{equation}
By Cardano's formula for the solutions of degree three equations we find
\begin{eqnarray*}
x & = & -\frac{\sqrt[3]{2802410\sqrt{5}+6281714+186\sqrt{2522364702+1128044106\sqrt{5}}}}{186}\\[3mm]
&&+\frac{2066/31+2986\sqrt{5}/93}{\sqrt[3]{2802410\sqrt{5}+6281714+186\sqrt{2522364702+1128044106\sqrt{5}}}}+\frac{163}{186}+\frac{53\sqrt{5}}{186} \\[3mm]
& \approx & 0,3944605378\quad\text{and}\\[2mm]
y & \approx & 0,2604339436\,.
\end{eqnarray*}

Similarly to the Cubus simus, equation~\eqref{EqXForDodSimum} implies that the vertices of the Dodecaedron simum cannot be constructed with ruler and compass on the face of the dodecahedron. Namely \eqref{EqXForDodSimum} is irreducible over the field $\BQ(\Phi)$. This can be seen in many ways. For example by reducing \eqref{EqXForDodSimum}, which lies in $\BZ[\Phi][x]$ modulo $2$ to $(\bar\Phi+1)x^3+\bar\Phi x^2+x+\bar\Phi\in\BF_2(\bar\Phi)[x]$. By plugging in all four elements of $\BF_2(\bar\Phi)=\BF_4$, one easily checks that the latter equation has no solution in $\BF_2(\bar\Phi)$, hence is irreducible over $\BF_2(\bar\Phi)$. Therefore also \eqref{EqXForDodSimum} is irreducible over $\BQ(\Phi)$ by \cite[Chapter~11, Theorems~3.9 and 4.2]{Artin} and $\BQ(\Phi)(x)$ is a field extension of degree $3$ over $\BQ(\Phi)$. This proves that $x$ cannot be constructed by ruler and compass but that it can be constructed by paper folding. 

Unfortunately we could not find a folding construction similarly elegant to our construction of the Cubus simus in Section~\ref{SectCubusSimus}. The problem is actually that there are infinitely many folding constructions and one wants to find the most elegant among them. In all constructions we found, the number $31$ plays an important role. We cannot explain this phenomenon but we observe that multiplying \eqref{EqXForDodSimum} with its conjugate (in which $\Phi=\frac{1+\sqrt{5}}{2}$ is replaced by $1-\Phi=\frac{1-\sqrt{5}}{2}$) and dividing by 31 yields
\[
\TS x^6-\frac{163}{31}\,x^5+\frac{386}{31}\,x^4-\frac{306}{31}\,x^3+\frac{89}{31}\,x^2-\frac{5}{31}\,x-\frac{1}{31}=0\,.
\]
This is irreducible over $\BQ$, because its reduction modulo $2$ is irreducible over $\BF_2$, and hence equals the minimal polynomial of $x$ over $\BQ$. Therefore the algebraic number $x\in\ol\BQ$ is integral at all primes different from $31$ but is not integral at $31$.

\begin{figure}[h]
\centering
\caption{Constructing the vertices of the Dodecaedron simum by paper folding.}
\label{FigureFoldingDodecaedrumSimum}
\psset{unit=0.8cm}
\begin{pspicture}(-9.5,-15)(8,10.5)
\psline[linewidth=0.5mm](0,0)(-4,0)
\psline[linewidth=0.5mm](-4,0)(-5.236,-3.804)
\psline[linewidth=0.5mm](-5.236,-3.804)(-2,-6.155)
\psline[linewidth=0.5mm](-2,-6.155)(1.236,-3.804)
\psline[linewidth=0.5mm](1.236,-3.804)(0,0)
\put(-0.4,0.1){$O$}
\put(-4.5,-0.1){$P$}
\put(-5.65,-3.9){$T$}
\put(-2.5,-6.5){$S$}
\put(1.3,-3.9){$R$}
\psline[linewidth=0.1mm](-4,0)(1.236,-3.804)
\psline[linewidth=0.1mm](5.767,9.401)(-2.078,-1.397)
\put(3.3,5.8){$g_1$}
\psdots[dotsize=1.5mm](-7.284,-7.314)(4.206,7.252)
\put(-7.1,-7.5){$G_1$}
\put(4.35,7.1){$G'_1$}
\psline[linestyle=dashed, linewidth=0.1mm](-7.284,-7.314)(4.206,7.252)
\psline[linewidth=0.1mm](5.767,0)(0,0)
\psline[linewidth=0.1mm](5.767,9.401)(5.767,0)
\put(5.9,2){$g_2$}
\psdots[dotsize=1.5mm](-7.409,-8.949)(5.767,7.755)
\put(-7.3,-9.3){$G_2$}
\put(5.9,7.6){$G'_2$}
\psline[linestyle=dashed, linewidth=0.1mm](-7.409,-8.949)(5.767,7.755)
\psline[linewidth=0.1mm](-2,-6.155)(-5.830,5.632)
\psline[linewidth=0.1mm](5.767,9.401)(-5.830,5.632)
\put(-1.3,7.4){$g_3$}
\psdots[dotsize=1.5mm](-8.801,-6.689)(3.246,8.582)
\put(-8.6,-6.8){$G_3$}
\put(3.1,8.1){$G'_3$}
\psline[linestyle=dashed, linewidth=0.1mm](-8.801,-6.689)(3.246,8.582)
\psline[linewidth=0.1mm](-2,-6.155)(3.314,10.198)
\psline[linewidth=0.1mm](5.767,9.401)(3.314,10.198)
\put(3.4,9.8){$g_4$}
\psdots[dotsize=1.5mm](-9.863,-7.938)(4.213,9.906)
\put(-9.7,-8.2){$G_4$}
\put(4.4,9.9){$G'_4$}
\psline[linestyle=dashed, linewidth=0.1mm](-9.863,-7.938)(4.213,9.906)
\psline[linewidth=0.1mm](-5.236,-3.804)(8.245,5.991)
\psline[linewidth=0.1mm](5.767,9.401)(8.245,5.991)
\put(7.2,7.8){$g_5$}
\psdots[dotsize=1.5mm](-9.003,-9.335)(5.772,9.395)
\put(-8.85,-9.6){$G_5$}
\put(5.95,9.4){$G'_5$}
\psline[linestyle=dashed, linewidth=0.1mm](-9.003,-9.335)(5.772,9.395)
\psline[linewidth=0.8mm](-7.481,4.656)(-6.602,3.963)(2.190,-2.972)(2.630,-3.319)
\put(2.3,-2.8){$f$}
\psline[linewidth=0.1mm](0,0)(-2,-6.155)
\psline[linestyle=dotted](-1.578,0)(-1.256,-0.991)
\psline[linestyle=dotted](-1.256,-0.991)(0.280,-0.991)
\psdot(-1.256,-0.991)
\pscircle(-0.8,-0.3){1.1}
\put(7,-8.5){
\pscircle(-4,-1.5){5.5}
\psdots[dotsize=2mm](0,0)(-7.889,0)(-1.610,-4.954)(-6.280,-4.954)
\put(0,0.2){$O$}
\put(-7.8,0.2){$X$}
\put(-5.8,0.2){$x$}
\put(-1.4,-4.8){$Y$}
\put(-1.2,-2.7){$y$}
\put(-6.8,-5.4){$A_1$}
\psline[linewidth=0.5mm,arrows=->,arrowsize=2mm,arrowinset=0](0,0)(-9.292,0)
\put(-9.3,-0.6){$P$}
\psline[linewidth=0.5mm,arrows=->,arrowsize=2mm,arrowinset=0](0,0)(1.131,-3.481)
\put(0.5,-3.7){$R$}
\psline[linewidth=0.8mm](-0.554,-5.787)(-8.971,0.853)
\put(-5.5,-2.6){$f$}
\psline[linewidth=0.1mm,arrows=->,arrowsize=2mm,arrowinset=0](0,0)(-2.173,-6.688)
\put(-2.6,-6.6){$S$}
\psline[linewidth=0.8mm,linestyle=dashed](-3.945,4)(-3.945,-7)
\put(-4.3,2){$\ell$}
\psline[linestyle=dotted](-7.889,0)(-6.280,-4.954)
\put(-7.4,-2.7){$y$}
\psline[linestyle=dotted](-6.280,-4.954)(0.280,-4.954)
\put(-5.3,-4.8){$x$}
}
\end{pspicture}
\end{figure}
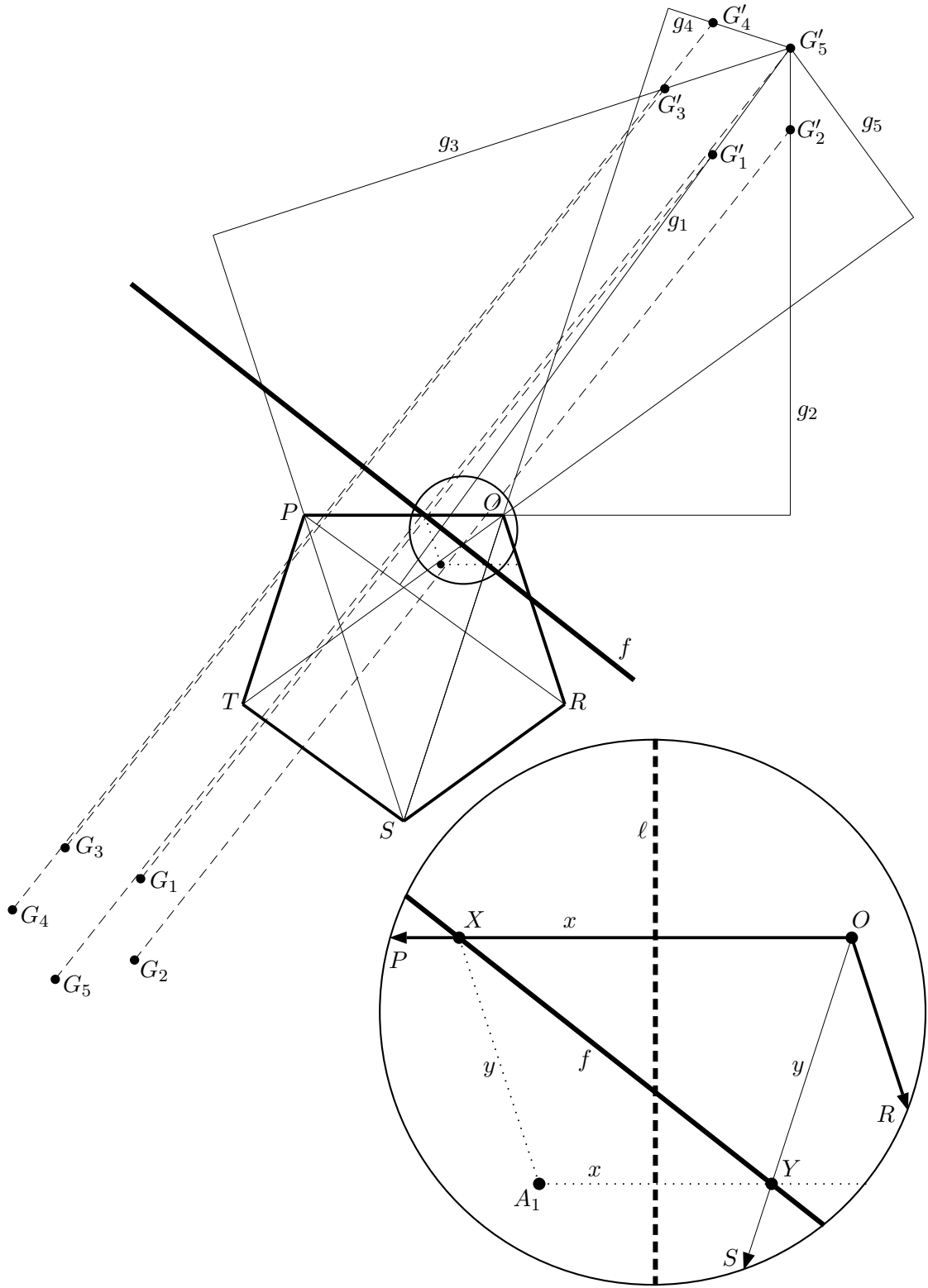

The nicest construction we found proceeds as follows. It creates a crease through the following two points $X$ and $Y$; see the magnifying glass in figure~\ref{FigureFoldingDodecaedrumSimum}. The point $X$ lies on the segment $[OP]$ at distance $x$ from $O$. The point $Y$ lies on the segment $[OS]$ at distance $y$ from $O$. This crease is created by placing (any two of) the points $G_1,G_2,G_3,G_4,G_5$ onto their respective line $g_1,g_2,g_3,g_4,g_5$ as in figure~\ref{FigureFoldingDodecaedrumSimum}.

\clearpage

In this figure the crease is the line $f$ and reflection at $f$ maps $G_i$ to the point $G'_i$ on the line $g_i$. The points $G_1,\ldots,G_5$ can be found as follows
\begin{eqnarray*}
\overrightarrow{OG_1}&=&\TS(\frac{37}{31}\Phi+\frac{15}{31})\cdot\overrightarrow{OP}+(\frac{22}{31}\Phi+\frac{24}{31})\cdot\overrightarrow{OR}\,,\\[2mm]
\overrightarrow{OG_2}&=&\TS(\frac{42}{31}\Phi+\frac{12}{31})\cdot\overrightarrow{OP}+(\frac{29}{31}\Phi+\frac{26}{31})\cdot\overrightarrow{OR}\,,\\[2mm]
\overrightarrow{OG_3}&=&\TS(\frac{47}{31}\Phi+\frac{9}{31})\cdot\overrightarrow{OP}+(\frac{17}{31}\Phi+\frac{27}{31})\cdot\overrightarrow{OR}\,,\\[2mm]
\overrightarrow{OG_4}&=&\TS(\frac{46}{31}\Phi+\frac{22}{31})\cdot\overrightarrow{OP}+(\frac{27}{31}\Phi+\frac{21}{31})\cdot\overrightarrow{OR}\,,\\[2mm]
\overrightarrow{OG_5}&=&\TS(\frac{49}{31}\Phi+\frac{14}{31})\cdot\overrightarrow{OP}+(\frac{26}{31}\Phi+\frac{34}{31})\cdot\overrightarrow{OR}\,.
\end{eqnarray*}
The line $g_1$ is the perpendicular to the line $PR$ at distance $\frac{12}{31}\Phi-\frac{1}{31}$ from $P$.\\[2mm]
The line $g_2$ is the perpendicular to the line $PO$ at distance $\frac{33}{62}\Phi+\frac{18}{31}$ from $O$.\\[2mm]
The line $g_3$ is the perpendicular to the line $SP$ at distance $\frac{11}{62}\Phi+\frac{37}{31}$ from $P$.\\[2mm]
The line $g_4$ is the perpendicular to the line $SO$ at distance $\frac{39}{31}\Phi+\frac{20}{31}$ from $O$.\\[2mm]
The line $g_5$ is the perpendicular to the line $TO$ at distance $\frac{38}{31}\Phi+\frac{35}{62}$ from $O$.\\[2mm]
Note that the points $G_1,\ldots,G_5$, and hence also the points $G'_1,\ldots,G'_5$ form regular pentagons. Also note that all the lines $g_1,\ldots,g_5$ have one point in common which is \emph{not} equal to $G'_5$ but very near so that the difference cannot be seen in figure~\ref{FigureFoldingDodecaedrumSimum}.

Once the crease $f$ has been constructed the point $A_1$ is obtained by one more folding along the dashed line $\ell$ (in the magnifying glass in figure~\ref{FigureFoldingDodecaedrumSimum}) which places the point $X$ onto $O$. This folding also places the point $Y$ onto $A_1$.

We found this (and many other) folding constructions using the computer algebra system Maple. We leave it to the interested reader to verify our construction, or even better, to find more elegant constructions for the Dodecaedron simum. We are happy to send our Maple code on request.

\begin{appendix}
\section{Verifying the construction rules from figures~\ref{FigureOctahedron}, \ref{FigureIcosahedron}, and \ref{FigureDodecahedron}.} \label{AppConstructions}

First note that the construction rules for the \emph{truncated tetrahedron}, the \emph{truncated cube}, the \emph{truncated octahedron}, the \emph{truncated icosahedron}, and the \emph{truncated dodecahedron}, as well as for the \emph{cuboctahedron}, and the \emph{icosidodecahedron} are well known; see for example~\cite{AdamWyss}. Likewise the rules to construct the \emph{vertices of the rhombicuboctahedron} and the \emph{truncated cuboctahedron on the face of a cube} are straight forward to prove.

\bigskip

\noindent\begin{minipage}[b]{0.67\linewidth}
\mbox{\hspace{5mm}} To verify the rule for the \emph{vertices of the rhombicuboctahedron on the face of an octahedron} and the \emph{vertices of the rhombicosidodecahedron on the face of an icosahedron} consider a triangle whose edge is divided into three parts with ratios $a:b:a$. Note that $|A''A'| = |A A'|$ as sides of the small rhombus. So we obtain for the ratio $\frac{|A'' B''|}{|A A'|} = \frac{|A' B'| - 2 |A'' A'|}{|A A'|} = \frac{|CA'| - 2 |AA'|}{|AA'|} = \frac{a+b}{a}-2=\frac{b-a}{a}$. 

\mbox{\hspace{5mm}} If the triangle is the face of an octahedron and $A''$ and $B''$ are the vertices of a rhombicuboctahedron, they form part of an octagon. In order that this octagon is regular we need that $\frac{|A''B''|}{|A''A'|}=\sqrt{2}$ (compare the picture for the truncated cube in figure~\ref{FigureCube}). This implies $b=(\sqrt{2}+1)\cdot a$ as desired.
\end{minipage}
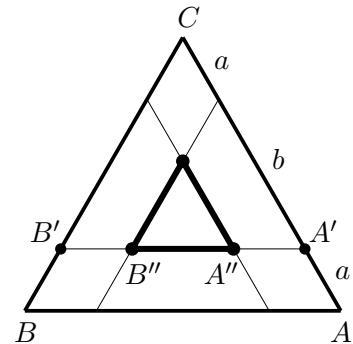
\begin{minipage}[b]{0.33\linewidth}
\psset{unit=0.7cm}
\begin{pspicture}(-1,6)(6,-0.5)
\psline[linewidth=0.5mm](3,5.196)(6,0)
\psline[linewidth=0.5mm](0,0)(3,5.196)
\psline[linewidth=0.5mm](0,0)(6,0)

\psline[linewidth=0.1mm](1.359,0)(3.680,4.019)
\psline[linewidth=0.1mm](2.320,4.019)(4.641,0)
\psline[linewidth=0.1mm](0.680,1.177)(5.320,1.177)

\psline[linewidth=0.8mm](2.039,1.177)(3.961,1.177)
\psline[linewidth=0.8mm](3.961,1.177)(3,2.842)
\psline[linewidth=0.8mm](3,2.842)(2.039,1.177)

\psdots[dotsize=1.8mm](2.039,1.177)(3.961,1.177)(3,2.842)

\psdots[dotsize=1.5mm](0.680,1.177)(5.320,1.177)
\put(5.8,-0.6){$A$}
\put(-0.2,-0.6){$B$}
\put(2.9,5.4){$C$}
\put(5.4,1.3){$A'$}
\put(0.1,1.3){$B'$}
\put(3.4,0.5){$A''$}
\put(1.9,0.5){$B''$}

\put(5.9,0.6){$a$}
\put(4.7,2.7){$b$}
\put(3.6,4.6){$a$}
\end{pspicture}
\end{minipage}

If the triangle is the face of an icosahedron and $A''$ and $B''$ are the vertices of a rhombicosidodecahedron, they form part of a decagon. In order that this decagon is regular we need that $|A''B''|=\Phi\cdot|A''A'|$ (compare the picture for the truncated dodecahedron in figure~\ref{FigureDodecahedron}). This implies $b=(\Phi+1)\cdot a$ as desired.

\vspace{-8mm}
\noindent
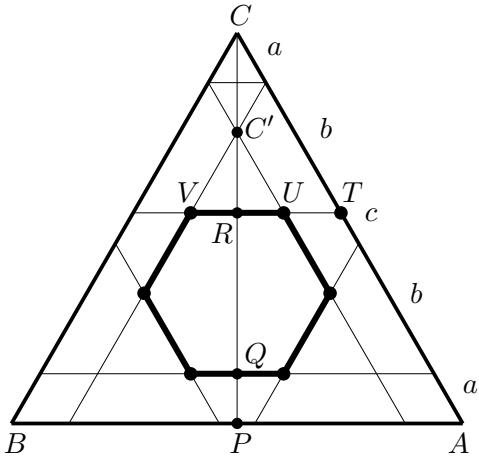
\begin{minipage}{0.4\linewidth}
\psset{unit=1cm}
\begin{pspicture}(0,7)(5.5,-1)
\psline[linewidth=0.5mm](3,5.196)(6,0)
\psline[linewidth=0.5mm](0,0)(3,5.196)
\psline[linewidth=0.5mm](0,0)(6,0)

\psline[linewidth=0.1mm](0.764,0)(3.382,4.535)
\psline[linewidth=0.1mm](2.618,4.535)(5.236,0)
\psline[linewidth=0.1mm](0.382,0.662)(5.618,0.662)

\psline[linewidth=0.1mm](3.236,0)(4.618,2.394)
\psline[linewidth=0.1mm](1.382,2.394)(2.764,0)
\psline[linewidth=0.1mm](1.618,2.803)(4.382,2.803)

\psline[linewidth=0.8mm](2.382,0.662)(3.618,0.662)
\psline[linewidth=0.8mm](3.618,0.662)(4.236,1.732)
\psline[linewidth=0.8mm](4.236,1.732)(3.618,2.803)
\psline[linewidth=0.8mm](3.618,2.803)(2.382,2.803)
\psline[linewidth=0.8mm](2.382,2.803)(1.764,1.732)
\psline[linewidth=0.8mm](1.764,1.732)(2.382,0.662)

\psdots[dotsize=1.8mm](2.382,0.662)(3.618,0.662)(4.236,1.732)(3.618,2.803)(2.382,2.803)(1.764,1.732)(4.382,2.803)

\psline[linewidth=0.1mm](3,0)(3,5.196)
\psline[linewidth=0.1mm](3.382,4.535)(2.618,4.535)

\psdots[dotsize=1.5mm](3,0.662)(3,2.803)(3,3.873)(3,0)
\put(5.8,-0.4){$A$}
\put(-0.1,-0.4){$B$}
\put(2.9,5.3){$C$}
\put(3.1,3.8){$C'$}
\put(2.9,-0.4){$P$}
\put(3.1,0.8){$Q$}
\put(2.65,2.4){$R$}
\put(4.4,2.95){$T$}
\put(3.6,2.95){$U$}
\put(2.2,2.95){$V$}

\put(6,0.4){$a$}
\put(5.3,1.6){$b$}
\put(4.7,2.7){$c$}
\put(4.1,3.8){$b$}
\put(3.4,4.9){$a$}
\end{pspicture}
\end{minipage}
\begin{minipage}{0.6\linewidth}
To construct the \emph{truncated cuboctahedron on the octahedron} and the \emph{truncated icosidodecahedron on the icosahedron} consider a triangle whose edge is divided into five parts with ratios $a\!:\!b\!:\!c\!:\!b\!:\!a$. We compute $\frac{|CC'|}{|CP|}=\frac{2a}{2a+2b+c}$, as well as $\frac{|CR|}{|CP|}=\frac{a+b}{2a+2b+c}$, and $\frac{|RQ|}{|CP|}=\frac{b+c}{2a+2b+c}$. This yields $|C'R|=|CR|-|CC'|=\frac{b-a}{2a+2b+c}\cdot|CP|$. Therefore the hexagon in the middle is regular if and only if $|RQ|=2|C'R|$, that is if and only if $c=b-2a$.

Furthermore, we compute $|UV|=|C'U|=|C'R|\cdot\frac{|AC|}{|CP|}=\frac{b-a}{2a+2b+c}\cdot|AC|$ and $|TU|=\frac{a}{2a+2b+c}\cdot|AB|$.
\end{minipage}

If the triangle is the face of an octahedron and $U$ and $V$ are the vertices of a truncated cuboctahedron, they form part of an octagon. In order that this octagon is regular we need that $|UV|=\sqrt{2}\cdot|TU|$ (compare the picture for the truncated cube in figure~\ref{FigureCube}). This implies $b=(\sqrt{2}+1)\cdot a$ and $c=(\sqrt{2}-1)\cdot a$ as desired.

If the triangle is the face of an icosahedron and $U$ and $V$ are the vertices of a truncated icosidodecahedron, they form part of a decagon. In order that this decagon is regular we need that $|UV|=\Phi\cdot|TU|$ (compare the picture for the truncated dodecahedron in figure~\ref{FigureDodecahedron}). This implies $b=(\Phi+1)\cdot a$ and $c=(\Phi-1)\cdot a$ as desired.

\medskip

To construct the vertices of a \emph{rhombicosidodecahedron on the face of a dodecahedron} divide the edges of a dodecahedron in three parts with ratios $1\!:\!1\!:\!1$.

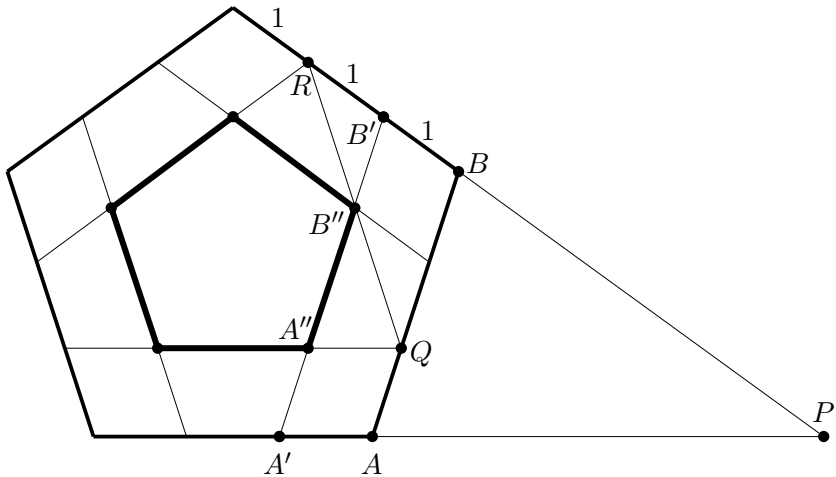
\psset{unit=1cm}
\begin{pspicture}(0,6.5)(11,-1)
\psline[linewidth=0.5mm](1.146,0)(4.854,0)
\psline[linewidth=0.5mm](4.854,0)(6,3.527)
\psline[linewidth=0.5mm](6,3.527)(3,5.706)
\psline[linewidth=0.5mm](3,5.706)(0,3.527)
\psline[linewidth=0.5mm](0,3.527)(1.146,0)

\psline[linewidth=0.1mm](4.854,0)(10.854,0)(6,3.527)
\psline[linewidth=0.1mm](4,4.980)(5.236,1.176)

\psline[linewidth=0.1mm](2.382,0)(1,4.253)
\psline[linewidth=0.1mm](3.618,0)(5,4.253)
\psline[linewidth=0.1mm](0.382,2.315)(4,4.980)
\psline[linewidth=0.1mm](0.764,1.176)(5.236,1.176)
\psline[linewidth=0.1mm](5.618,2.315)(2,4.980)

\psline[linewidth=0.8mm](2,1.176)(4,1.176)
\psline[linewidth=0.8mm](4,1.176)(4.618,3.041)
\psline[linewidth=0.8mm](4.618,3.041)(3,4.254)
\psline[linewidth=0.8mm](3,4.254)(1.382,3.041)
\psline[linewidth=0.8mm](1.382,3.041)(2,1.176)

\psdots[dotsize=1.5mm](1.382,3.041)(2,1.176)(4,1.176)(4.618,3.041)(3,4.254)
\psdots[dotsize=1.5mm](4.854,0)(10.854,0)(6,3.527)(3.618,0)(5,4.253)(4,4.980)(5.236,1.176)(3.618,0)(5,4.253)
\put(3.5,5.45){$1$}
\put(4.5,4.7){$1$}
\put(5.5,3.95){$1$}

\put(4.7,-0.5){$A$}
\put(6.1,3.5){$B$}
\put(3.4,-0.5){$A'$}
\put(4.5,3.9){$B'$}
\put(3.6,1.3){$A''$}
\put(4.0,2.7){$B''$}
\put(10.7,0.2){$P$}
\put(5.35,1.0){$Q$}
\put(3.75,4.55){$R$}
\end{pspicture}

We compute $|PB|=\Phi\cdot|AB|$ and $|A'B'|=\frac{1/3+\Phi}{\Phi}\cdot|AB|=(1+\frac{1}{3\Phi})\cdot3|BB'|=(2+\Phi)|BB'|$. Since $|B'B''|=|BB'|$ as sides of a rhombus, $|A''B''|=|A'B'|-2|B'B''|=\Phi\cdot|BB'|$. Moreover, $|QR|=\Phi|BR|=2\Phi\cdot|BB'|$. The point $B''$ together with two more vertices of the rhombicosidodecahedron forms an equilateral triangle of edge length $\frac{1}{2}|QR|=\Phi|BB'|=|A''B''|$. So all edges of the rhombicosidodecahedron have the same length and this proves that our construction is correct.


\bigskip

Finally look at the face of a dodecahedron and divide its edges into five parts with lengths $a,b,c,b,a$ to construct the vertices of the \emph{truncated icosidodecahedron}.

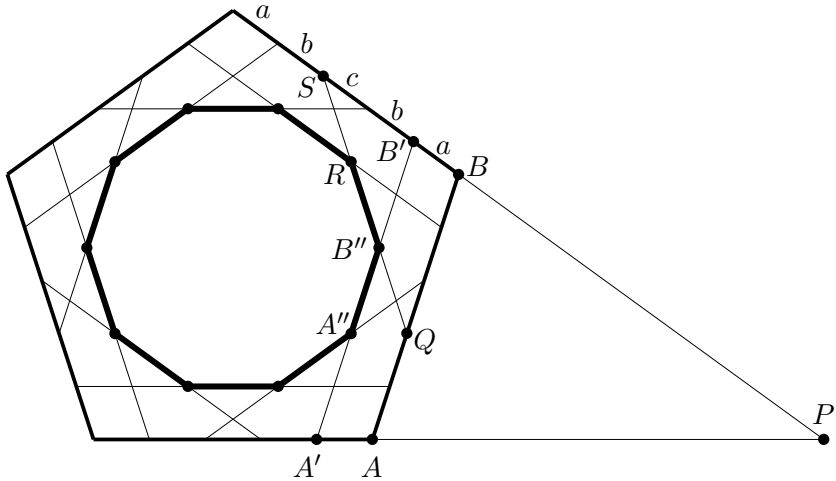
\psset{unit=1cm}
\begin{pspicture}(0,6.5)(11,-1)
\psline[linewidth=0.5mm](1.146,0)(4.854,0)
\psline[linewidth=0.5mm](4.854,0)(6,3.527)
\psline[linewidth=0.5mm](6,3.527)(3,5.706)
\psline[linewidth=0.5mm](3,5.706)(0,3.527)
\psline[linewidth=0.5mm](0,3.527)(1.146,0)

\psline[linewidth=0.1mm](4.854,0)(10.854,0)(6,3.527)

\psline[linewidth=0.1mm](1.888,0)(0.6,3.963)
\psline[linewidth=0.1mm](4.112,0)(5.4,3.963)
\psline[linewidth=0.1mm](0.229,2.822)(3.6,5.270)
\psline[linewidth=0.1mm](0.917,0.705)(5.083,0.705)
\psline[linewidth=0.1mm](5.771,2.822)(2.4,5.270)

\psline[linewidth=0.1mm](2.629,0)(5.542,2.116)
\psline[linewidth=0.1mm](3.371,0)(0.458,2.116)
\psline[linewidth=0.1mm](0.688,1.411)(1.8,4.834)
\psline[linewidth=0.1mm](1.2,4.399)(4.8,4.399)
\psline[linewidth=0.1mm](5.312,1.411)(4.2,4.834)

\psline[linewidth=0.8mm](2.4,0.705)(3.6,0.705)
\psline[linewidth=0.8mm](3.6,0.705)(4.570,1.410)
\psline[linewidth=0.8mm](4.570,1.410)(4.942,2.553)
\psline[linewidth=0.8mm](4.942,2.553)(4.571,3.693)
\psline[linewidth=0.8mm](4.571,3.693)(3.6,4.399)
\psline[linewidth=0.8mm](3.6,4.399)(2.4,4.399)
\psline[linewidth=0.8mm](2.4,4.399)(1.429,3.693)
\psline[linewidth=0.8mm](1.429,3.693)(1.058,2.553)
\psline[linewidth=0.8mm](1.058,2.553)(1.430,1.410)
\psline[linewidth=0.8mm](1.430,1.410)(2.4,0.705)

\psdots[dotsize=1.5mm](2.4,0.705)(3.6,0.705)(4.570,1.410)(4.942,2.553)(4.571,3.693)(2.4,4.399)(3.6,4.399)(1.429,3.693)(1.058,2.553)(1.430,1.410)
\psdots[dotsize=1.5mm](4.854,0)(6,3.527)(4.112,0)(5.4,3.963)(10.854,0)(5.312,1.411)(4.2,4.834)
\put(3.3,5.6){$a$}
\put(3.9,5.15){$b$}
\put(4.5,4.7){$c$}
\put(5.1,4.25){$b$}
\put(5.7,3.8){$a$}

\put(4.7,-0.5){$A$}
\put(6.1,3.5){$B$}
\put(3.8,-0.5){$A'$}
\put(4.9,3.7){$B'$}
\put(4.1,1.4){$A''$}
\put(4.3,2.4){$B''$}
\put(10.7,0.2){$P$}
\put(5.4,1.2){$Q$}
\put(4.2,3.4){$R$}
\put(3.85,4.55){$S$}
\end{pspicture}

This time $|PB|=\Phi\cdot|AB|$ and $|A'B'|=\frac{a+\Phi\cdot|AB|}{\Phi\cdot|AB|}\cdot|AB|=|AB|+\frac{a}{\Phi}=2a+2b+c+a\Phi-a=a+2b+c+a\Phi$. Also $|B'B''|=b+c$ as sides of an isosceles triangle. Therefore $|A''B''|=|A'B'|-2|B'B''|=a-c+a\Phi$. On the other hand $|QS|=\Phi\cdot(a+b+c)$ and $|B''Q|=\frac{a}{a+b+c}\cdot|QS|=a\Phi$ imply $|B''R|=|QS|-2|B''Q|=\Phi\cdot(b+c-a)$. The points $B''$ and $R$ together with four other vertices of the truncated icosidodecahedron form a hexagon. It is regular if and only if $|B''R|=|B''Q|$ (compare the picture for the truncated icosahedron in figure~\ref{FigureIcosahedron}), that is, if and only if $b+c=2a$. In addition, the decagon containing $A'',B''$ and $R$ is regular if and only if $|A''B''|=|B''R|=|B''Q|$, which is equivalent to $a=c$. Thus we find $a=c=b$. 

\end{appendix}

%
%

\vfill

{\small
\begin{minipage}[t]{0.5\linewidth}
\noindent
Urs Hartl\\
Universit\"at M\"unster\\
Mathematisches Institut \\
Einsteinstr.~62\\
D -- 48149 M\"unster
\\ Germany
\\[1mm]
\href{http://www.math.uni-muenster.de/u/urs.hartl/}{\small www.math.uni-muenster.de/u/urs.hartl/}
\end{minipage}
\begin{minipage}[t]{0.45\linewidth}
\noindent
Klaudia Kwickert\\
Universit\"at M\"unster\\
Mathematisches Institut \\
Einsteinstr.~62\\
D -- 48149 M\"unster
\\ Germany
\end{minipage}
}

\end{document}